\documentclass[11pt]{article}
\usepackage{cancel}
\usepackage{hyperref}
\usepackage{amssymb}
\usepackage{amsmath}
\usepackage{latexsym}

\usepackage{amsthm}

\font\tengoth=eufm10 at 10pt
\font\sevengoth=eufm7 at 6pt
\newfam\gothfam
\textfont\gothfam=\tengoth
\scriptfont\gothfam=\sevengoth

\newcommand{\mlabel}[1]{\marginpar{#1}\label{#1}}

\newcommand{\g}{{\mathfrak g}}

\newcommand{\fa}{{\mathfrak a}}
\newcommand{\fb}{{\mathfrak b}}
\newcommand{\fc}{{\mathfrak c}}

\newcommand{\fe}{{\mathfrak e}}

\newcommand{\fg}{{\mathfrak g}}
\newcommand{\fh}{{\mathfrak h}}

\newcommand{\fk}{{\mathfrak k}}
\newcommand{\fl}{{\mathfrak l}}

\newcommand{\fn}{{\mathfrak n}}

\newcommand{\fq}{{\mathfrak q}}
\newcommand{\fp}{{\mathfrak p}}
\newcommand{\fr}{{\mathfrak r}}
\newcommand{\fs}{{\mathfrak s}}
\newcommand{\ft}{{\mathfrak t}}
\newcommand{\fu}{{\mathfrak u}}

\newcommand{\fz}{{\mathfrak z}}

\renewcommand\sp{\mathfrak {sp}} 
 
\newcommand\conf{\mathfrak {conf}}

\renewcommand{\:}{\colon}
\newcommand{\1}{\mathbf{1}}

\newcommand{\cE}{\mathcal{E}}

\newcommand{\cG}{\mathcal{G}}
\newcommand{\cH}{\mathcal{H}}

\newcommand{\cM}{\mathcal{M}}

\newcommand{\cO}{\mathcal{O}}

\newcommand{\cW}{\mathcal{W}}

\newcommand{\bO}{\mathbf{O}}

\newcommand{\eset}{\emptyset}

\newcommand{\dd}{{\tt d}}

\newcommand{\trile}{\trianglelefteq}
\newcommand{\subeq}{\subseteq}
\newcommand{\supeq}{\supseteq}

\newcommand{\into}{\hookrightarrow}
\newcommand{\eps}{\varepsilon}

\newcommand{\shalf}{{\textstyle{\frac{1}{2}}}}

\def\onto{\to\mskip-14mu\to}

\newcommand{\N}{{\mathbb N}}
\newcommand{\Z}{{\mathbb Z}}
\newcommand{\R}{{\mathbb R}}
\newcommand{\C}{{\mathbb C}}

\newcommand{\K}{{\mathbb K}}

\renewcommand{\H}{{\mathbb H}}

\newcommand{\bH}{{\mathbb H}}
\newcommand{\bS}{{\mathbb S}}

\renewcommand{\tilde}{\widetilde}

\newcommand{\GL}{\mathop{{\rm GL}}\nolimits}
\newcommand{\SL}{\mathop{{\rm SL}}\nolimits}

\newcommand{\PSL}{\mathop{{\rm PSL}}\nolimits}
\newcommand{\SO}{\mathop{{\rm SO}}\nolimits}
\newcommand{\SU}{\mathop{{\rm SU}}\nolimits}
\newcommand{\OO}{\mathop{\rm O{}}\nolimits}

\newcommand{\U}{\mathop{\rm U{}}\nolimits}

\newcommand{\Sym}{\mathop{{\rm Sym}}\nolimits}

\newcommand{\Skew}{\mathop{{\rm Skew}}\nolimits}

\newcommand{\gl}  {\mathop{{\mathfrak{gl} }}\nolimits}

\newcommand{\fsl} {\mathop{{\mathfrak{sl} }}\nolimits}
\newcommand{\fsp} {\mathop{{\mathfrak{sp} }}\nolimits}
\newcommand{\su}  {\mathop{{\mathfrak{su} }}\nolimits}
\newcommand{\so}  {\mathop{{\mathfrak{so} }}\nolimits}

\newcommand{\Fix}{\mathop{{\rm Fix}}\nolimits}

\newcommand{\ad}{\mathop{{\rm ad}}\nolimits}
\newcommand{\Ad}{\mathop{{\rm Ad}}\nolimits}

\newcommand{\tr}{\mathop{{\rm tr}}\nolimits}

\newcommand{\Herm}{\mathop{{\rm Herm}}\nolimits}
\newcommand{\Aherm}{\mathop{{\rm Aherm}}\nolimits}

\newcommand{\Aut}{\mathop{{\rm Aut}}\nolimits}

\newcommand{\Conf}{\mathop{\rm Conf{}}\nolimits}

\newcommand{\diag}{\mathop{{\rm diag}}\nolimits}

\newcommand{\id}{\mathop{{\rm id}}\nolimits}

\renewcommand{\dim}{\mathop{{\rm dim}}\nolimits}

\newcommand{\Inn}{\mathop{{\rm Inn}}\nolimits}

\newcommand{\Spann}{\mathop{{\rm span}}\nolimits}

\newcommand{\Str}{\mathop{{\rm Str}}\nolimits}
\newcommand{\dS}{\mathop{{\rm dS}}\nolimits}

\newcommand{\Rarrow}{\Rightarrow}
\newcommand{\nin}{\noindent} 
\newcommand{\oline}{\overline}

\newcommand{\la}{\langle}
\newcommand{\ra}{\rangle}

\newcommand{\res}{\vert}

\newcommand{\Spec}{{\rm Spec}}
\newcommand{\Spin}{{\rm Spin}}

\newcommand{\ssssarr}{\hbox to 15pt{\rightarrowfill}}
\newcommand{\sssarr}{\hbox to 20pt{\rightarrowfill}}
\newcommand{\ssarr}{\hbox to 30pt{\rightarrowfill}}
\newcommand{\sarr}{\hbox to 40pt{\rightarrowfill}}
\newcommand{\arr}{\hbox to 60pt{\rightarrowfill}}
\newcommand{\sssslarr}{\hbox to 15pt{\leftarrowfill}}
\newcommand{\ssslarr}{\hbox to 20pt{\leftarrowfill}}
\newcommand{\sslarr}{\hbox to 30pt{\leftarrowfill}}
\newcommand{\slarr}{\hbox to 40pt{\leftarrowfill}}
\newcommand{\larr}{\hbox to 60pt{\leftarrowfill}}

\newcommand{\Arr}{\hbox to 80pt{\rightarrowfill}}

\def\theoremname{Theorem}
\def\propositionname{Proposition}
\def\corollaryname{Corollary}
\def\lemmaname{Lemma}
\def\remarkname{Remark}
\def\conjecturename{Conjecture} 

\def\definitionname{Definition}
\def\exercisename{Exercise}
\def\examplename{Example}
\def\examplesname{Examples}
\def\problemname{Problem}
\def\problemsname{Problems}

\def\satzname{Satz} 
\def\koroname{Korollar}
\def\folgname{Folgerung}
\def\bemerkname{Bemerkung}
\def\aufgname{Aufgabe}

\def\beisname{Beispiel}
\def\beissname{Beispiele}
\def\bewname{Beweis}

\def\@thmcounter#1{\noexpand\arabic{#1}}
\def\@thmcountersep{}
\def\@begintheorem#1#2{\it \trivlist \item[\hskip 
\labelsep{\bf #1\ #2.\quad}]}
\def\@opargbegintheorem#1#2#3{\it \trivlist
      \item[\hskip \labelsep{\bf #1\ #2.\quad{\rm #3}}]}
\makeatother
\newtheorem{theor}{\theoremname}[section]
\newtheorem{propo}[theor]{\propositionname}
\newtheorem{coro}[theor]{\corollaryname}
\newtheorem{lemm}[theor]{\lemmaname}

\newenvironment{thm}{\begin{theor}\it}{\end{theor}}
\newenvironment{theorem}{\begin{theor}\it}{\end{theor}}

\newenvironment{prop}{\begin{propo}\it}{\end{propo}}

\newenvironment{cor}{\begin{coro}\it}{\end{coro}}

\newenvironment{lem}{\begin{lemm}\it}{\end{lemm}}

\newtheorem{rema}[theor]{\remarkname}

\newenvironment{remark}{\begin{rema}\rm}{\end{rema}}
\newenvironment{rem}{\begin{rema}\rm}{\end{rema}}

\newtheorem{stepnow}[theor]{}

\newtheorem{defin}[theor]{\definitionname} 

\newenvironment{defn}{\begin{defin}\rm}{\end{defin}}

\newtheorem{exerc}{\exercisename}[section]

\newtheorem{exa}[theor]{\examplename}

\newenvironment{ex}{\begin{exa}\rm}{\end{exa}}

\newtheorem{exas}[theor]{\examplesname}

\newtheorem{conj}[theor]{\conjecturename}

\newtheorem{pro}[theor]{\problemname}

\newtheorem{prs}[theor]{\problemsname}

\newtheorem{aufg}{\aufgname}[section]

\newenvironment{prf}{\begin{proof}}{\end{proof}}
 
\newcommand{\pmat}[1]{\begin{pmatrix} #1 \end{pmatrix}}

{\hfill\qed\end{trivlist}}

\newenvironment{beweis*}{\begin{trivlist}\item[\hskip%
\labelsep{\bf\bewname.\quad}]}%
{\end{trivlist}}

\newtheorem{satzn}[theor]{\satzname}

\newtheorem{koro}[theor]{\koroname}

\newtheorem{folg}[theor]{\folgname}

\newtheorem{bem}[theor]{\bemerkname}

\newtheorem{aufgn}[theor]{\aufgname}

\newtheorem{beis}[theor]{\beisname}

\newtheorem{beiss}[theor]{\beissname}

\numberwithin{equation}{section}

\usepackage{color}

\newcommand\be{{\bf{e}}}

\newcommand\RR{{\mathbb R}}

\newcommand{\detE}{{\textstyle{\det_E}}}

\newcommand{\sV}{{\tt V}}

\newcommand{\dom}{\mathop{{\rm dom}}\nolimits}

\newcommand{\oEp}{\mathcal{E}^{\bot}}
\newcommand{\SC}{\mathcal{SC}}
\newcommand{\POh}{\mathrm{P}\mathcal{O}_h}

\renewcommand{\bO}{\mathbb O}

\renewcommand{\phi}{\varphi}

\renewcommand\mlabel{\label}

\begin{document}

\title{Orthogonal pairs of Euler elements I. \\
  Classification, fundamental groups and twisted duality}
\author{Vincenzo Morinelli, Karl-Hermann Neeb, Gestur \'Olafsson} 

\maketitle

\abstract{The current article continues our project
  on representation theory, Euler elements, causal homogeneous spaces and
  Algebraic Quantum Field Theory (AQFT).
  We call a pair $(h,k)$ of Euler elements
  orthogonal if $e^{\pi i \ad h} k = -k$.
  We show that,   if $(h,k)$ and $(k,h)$ are orthogonal, then
  they generate a $3$-dimensional simple subalgebra. 
  We also classify orthogonal Euler pairs in simple Lie algebras
  and determine the fundamental groups of adjoint Euler elements
  in arbitrary finite-dimensional Lie algebras.
  Causal complements of wedge regions in spacetimes
  can be related to so-called
  twisted complements in the space
  of  abstract Euler wedges, defined in purely group theoretic terms.
  We show that any pair of twisted complements can be connected
  by a chain of successive complements coming from $3$-dimensional subalgebras. 
}

\tableofcontents 

\vspace{1cm}  
  
\section{Introduction}

The current article continues our project
  on representation theory, causal homogeneous spaces and
  Algebraic Quantum Field Theory (AQFT)
  \cite{MN21, NO23a, NO23b, MNO23, MNO24, Mo25}.
Here we deal with algebraic structures related to finite-dimensional
Lie groups that arise naturally in
Algebraic Quantum Field Theory (AQFT). Specifically, we consider
{\it Euler elements} $h$ in a Lie algebra $\g$, i.e., elements for which
$\ad h$ is non-zero and diagonalizable with eigenvalues contained in $\{-1,0,1\}$.
Then $\tau_h^\g := e^{\pi i \ad h}$ is an involutive automorphism of $\g$
and we consider a connected Lie group $G$ with Lie algebra $\g$ for which $\tau_h^\g$
integrates to an automorphism $\tau_h$ of $G$
(which is always the case if $G$ is simply connected). We then form
the extended group
\begin{equation}
  \label{eq:gtauh}
  G_{\tau_h} := G\rtimes \{ \1,\tau_h\}.
  \end{equation} 

  We now indicate how these structures appear in AQFT.
In AQFT, in the sense of Haag--Kastler, one considers 
families (also  called nets) of von Neumann algebras $\cM(\cO)$
on a fixed Hilbert space $\cH$,
associated to open subsets $\cO$ in some space-time manifold~$M$ 
(\cite{Ha96}, \cite{Ne25}). The hermitian elements of the algebra $\cM(\cO)$ represent
observables  that can be measured in the ``laboratory'' $\cO$.
In this context the technique of modular localization
requires a unitary representation $U \: G \to \U(\cH)$ of a
(connected) symmetry group $G$ and a $G$-fixed unit vector $\Omega$
(representing a vacuum state), which for some algebras
$\cM(W)$ is cyclic ($\cM(W)\Omega$ is dense) and
separating (the orbit map $A \mapsto A\Omega$ is injective on
$\cM(W)$).  Then the Tomita--Takesaki Theorem
(\cite[Thm.~2.5.14]{BR87}) implies
the existence of a positive operator $\Delta$ (called the modular operator)
and a conjugation $J$ such that
\begin{equation}
  \label{eq:tomtak}
  J \cM(W) J = \cM(W)'  \quad \mbox{ and } \quad
 \Delta^{it} \cM(W) \Delta^{-it} = \cM(W) \quad \mbox{ for } \quad t \in \R. 
\end{equation}
The pairs $(\Delta, J)$ appearing in this context are characterized
by the modular relation
\begin{equation}
  \label{eq:modrel}
  J \Delta J = \Delta^{-1}.
\end{equation}
A key link between symmetries and the modular structures is implemented by
the Bisognano--Wichmann (BW) condition that requires the existence of a
Lie algebra element $h \in \g$, for which
the infinitesimal generator (a skew-adjoint operator) 
\[ \partial U(h) := \frac{d}{dt}\Big|_{t = 0}  U(\exp th) \]
represents the modular operator by
\[ \Delta = e^{2\pi i \cdot \partial U(h)}.\]
This implies in particular that $J \partial U(h) J  =\partial U(h)$ because
$J$ is antilinear. 
From \cite{MN24} we know that,
if $\ker U$ is discrete, this further implies that $h$ is an Euler element
and that, if $\tau_h$ exists on $G$, then $U$ extends by
$U(\tau_h) := J_\sV$ to an antiunitary representation of~$G_{\tau_h}$ on~$\cH$.

On the Lie group level, an analog of the pairs
$(\Delta, J)$ satisfying the modular relation \eqref{eq:modrel} 
is the {\it abstract wedge space}
  \begin{equation}
    \label{eq:cGGtau}
 \cG(G_{\tau_h}) := \{ (x,\sigma)\in \g \times G\tau_h
 \subeq \g \times G_{\tau_h}    \: \sigma^2 = e, \Ad(\sigma)x = x\}.
  \end{equation}
In view of \cite{MN24}, the most relevant ones are the {\it
  abstract Euler wedges}, 
denoted $\cG_E(G_{\tau_h})$, for which
$x$ is an Euler element and $\Ad(\sigma) = e^{\pi i \ad x} = \tau_x^\g$.
The set $\cG(G_{\tau_h})$ carries a natural 
$G$-action by $g.(x,\sigma)= (\Ad(g)x, g\sigma g^{-1})$
and a duality operation $(x,\sigma)' := (-x,\sigma)$, which
on the level of operator algebras corresponds to the passage to the
commutant (cf.\ \eqref{eq:tomtak} and \eqref{eq:modrel}).

The orbit
\begin{equation}
  \label{eq:cw+}
  \cW_+ := G.(h,\tau_h) \subeq \cG_E(G_{\tau_h})
\end{equation}
of the basic pair $W_0 := (h,\tau_h)$ thus encodes in an abstract
$G$-equivariant form the modular operators arising in
AQFTs. One therefore expects to understand 
inclusions and locality properties of the family of algebras $\cM(\cO)$ 
in terms of the geometry of the
abstract wedge space~$\cG_E(G_{\tau_h})$
(see \cite{BGL02}, \cite{NO17}, \cite{MN21} for more details and background).
If $Z(G)$ is discrete, then projecting onto the first factor yields a
$G$-equivariant covering map
\[ q_h^G \: \cW_+ \to \cO_h.\]
The fiber over $h$ is the orbit of $(h,\tau_h)$ under the stabilizer
$G^h$, whose elements are parametrized by the image
$Z_2 := \delta_h(G^h)$ of $G^h$ under the homomorphism
$\delta_h(g) = g \tau_h(g)^{-1}$, which plays the role of a connecting
  homomorphism in cohomology, hence the notation. 

To understand locality properties of $\cW_+$, we
observe that $h$ is {\it symmetric} in the sense that $-h \in
\cO_h := \Ad(G)h$ if and only if $\cW_+$ contains a {\it twisted complement} 
\[ W^{'\alpha} = (-h,\alpha \tau_h), \quad \alpha \in Z(G).\]
The possible complements of the abstract wedge
$(h,\tau_h) \in \cG_E(G_{\tau_h})$ are para\-met\-rized by the set
\[ \delta_h(\{  g \in G \: \Ad(g)h = -h\}) \subeq Z(G).\]
If  $G$ is simply connected, then $\cW_+$ is simply connected and 
\[ Z_2 \cong \pi_1(\cO_h) \] is the fundamental group
of the adjoint orbit $\cO_h$.
One of our central results asserts that,
if $\g$ is simple and $h \in \g$ is an
Euler elements, then $\pi_1(\cO_h)$ is either trivial, $\Z$ or $\Z_2$.
More specifically, infinite fundamental
groups occur precisely in the Cayley type case, i.e., when $\g$ is hermitian,
and $\Z_2$ occurs for the so-called split-type cases
(cf.\ \cite[Table 3]{MNO23}, Theorem~\ref{thm:Z2-struc}).
For general Lie algebras, the calculation of $\pi_1(\cO_h)$
can be reduced by Levi projections to the semisimple case,
which is a product of simple cases (Theorem~\ref{thm:leviproj}). \\

A key philosophy of modular theory in AQFT is to understand
symmetry groups which are generated by modular groups of von Neumann
algebras. This requires a good understanding of specific configurations
of such modular groups (\cite{Bo00}, \cite{Wi93, Wi97}, \cite{NO17}, \cite{Ne18}),
and such
configurations can be studied on the abstract level in $\cG_E(G_{\tau_h})$.
This led us to the subject matter of this papers, namely pairs
$(h, k)$ of Euler elements, which are {\it orthogonal}
in the sense that $\tau_h^\g(k) = -k$.\begin{footnote}
  {If $\kappa$ is an invariant symmetric bilinear form on $\g$, then
      it is invariant under $\tau_h^\g$, so that the relations
      $\tau_h^\g(h) = h$ and $\tau_h^\g(k) = -k$ imply that $h$ and $k$
      are $\kappa$-orthogonal. }  
\end{footnote}
A typical example of an 
  orthogonal pair of  Euler elements in $\fsl_2(\R)$ is $(h_0, k_0)$ with 
\begin{equation}
  \label{eq:h0k0-intro}
 h_0 := \frac{1}{2}\pmat{1 & 0 \\ 0 &-1} \quad \mbox{ and } \quad 
 k_0 := \frac{1}{2}\pmat{0 & 1 \\ 1 &0}.
\end{equation}
It is easy to see that any Euler element $h$ contained in some
$\fsl_2(\R)$-subalgebra $\fs \subeq \g$ 
is symmetric and admits a partner $k$ for which $(h,k)$ and $(k,h)$ are orthogonal.
Conversely, for any such pair $h$, $k$ and $z_{h,k} := [h,k]$ span a
subalgebra $\fs$ isomorphic to $\fsl_2(\R)$  (Theorem~\ref{thm:sl2-gen}). 
This underlines the prominent role of $\fsl_2(\R)$ in our context,
and can even be used to shed new light on twisted locality conditions.
For $h,k$ as above,
\begin{equation}
  \label{eq:zetahkintro}
  \zeta_{h,k} := \exp(2\pi z_{h,k}) \in Z(G)
\end{equation}
generates the cyclic center of the integral
subgroup  $S := \la \exp \fs \ra \subeq G$,
and $g = \exp(\pi z_{h,k}) \in S$ satisfies
\[ \Ad(g)h = -h \quad \mbox{ and }  \quad g.(h, \tau_h) = (-h, \zeta_{h,k} \tau_h).\]
We conclude that twisted complements
of abstract Euler wedge can be obtained from $\SL_2(\R)$-subgroups.

The main results of the present paper are the following: 
\begin{itemize}
\item A classification of all orthogonal pairs of Euler elements
  in simple Lie algebras under the action of the group
  $\Inn(\g) = \Ad(G)$ of inner automorphisms (Theorem~\ref{thm:1.4}). 
It is based on our classification of Euler elements 
in \cite{MN21} and structural results on $3$-graded
simple Lie algebras.
\item An extension of some of the results of \cite{MN21} from simple
  to general Lie algebras, which assert
  that an Euler element $h \in \g$ is symmetric if and only
  if there exists a second Euler element $k$ for which
  $(h,k)$ and $(k,h)$ are orthogonal. This in turn is equivalent to
  $h$ and $k$ generating a subalgebra isomorphic to $\fsl_2(\R)$
    (Theorem~\ref{thm:sl2-gen})
\item The calculation of the fundamental groups
    $\pi_1(\cO_h)$ for all (not necessarily symmetric)
    Euler elements by reduction to simple Lie algebras
    (Theorems~\ref{thm:Z2-struc} and \ref{thm:leviproj}). 
  \item All twisted complements $W_0^{'\alpha}$ can be obtained by
    successively acting with central elements of suitable $\SL_2(\R)$-subgroups.
    More specifically, the elements $\zeta_{h,k}$, for $h$ fixed and
    $(h,k), (k,h)$ orthogonal, generate the central subgroup
    $Z_3 := \delta_h(G^{\{\pm h\}})$  (Theorem~\ref{thm:slconj}). 
  \end{itemize}

  The structure of this paper is as follows. In Section~\ref{sec:2}
  we first discuss orthogonal Euler pairs and prove
  Theorem~\ref{thm:sl2-gen}, relating them to $\fsl_2$-subalgebras.
  The classification of orthogonal Euler pairs is given in
  Subsection~\ref{subsec:2.2}, and in Subsection~\ref{subsec:2.3}
  we characterize those that are symmetric in the sense that
  $(h,k)$ is conjugate to $(k,h)$. 
  The abstract wedge space and
  the central elements \eqref{eq:zetahkintro} are  introduced in
Section~\ref{sec:4}, 
  and in Section~\ref{sec:5} we calculate the fundamental groups $\pi_1(\cO_h)$.
  Twisted complements of abstract wedges are studied in Section~\ref{sec:6}.
  Here our main result is that all twists can be obtained by
  successively applying central elements $\zeta_{h,k}$
  contained in $\SL_2(\R)$-subgroups (Theorem~\ref{thm:slconj}).
The Outlook section connects the present paper with the
  forthcoming work in~\cite{MNO25}.
  We also added four appendices in which we provide some background
  on Jordan algebras, the relation between orthogonal
  Euler pairs and nilpotent elements,
 the conformal  Lie algebra $\so_{2,d}(\R)$ of Minkowski space $\R^{1,d-1}$,
    an alternative approach
  to $\pi_1(\cO_h)$ based on Wiggerman's presentation of
  the fundamental group of a flag manifold, and finally a
  relation with the kernel of the universal complexification
  $\eta_G \: G \to G_\C$ of a simply connected simple Lie group~$G$. \\  
  
  In a sequel \cite{MNO25} to this paper we describe
  the consequences of the algebraic results obtained in this paper
  for nets of local algebras and for the associated nets of real subspaces
  (cf.\ \cite{Ne25}).

\vspace{5mm}

\nin {\bf Notation and conventions:} 
\begin{itemize}
\item Throughout this paper $G$ denotes a finite-dimensional
  Lie group. 
\item We write $e \in G$ for the identity element in the Lie group~$G$ 
and $G_e$ for its identity component. 
\item For $x \in \g$, we write $G^x := \{ g \in G \: \Ad(g)x = x \}$ 
for the stabilizer of $x$ in the adjoint representation 
and $G^x_e = (G^x)_e$ for its identity component. 
\item For $h \in \g$ and $\lambda \in \R$, we write 
$\g_\lambda(h) := \ker(\ad h - \lambda \1)$ for the corresponding eigenspace 
in the adjoint representation.
\item For a Lie subalgebra $\fs \subeq \g$, we write 
$\Inn_\g(\fs)= \la e^{\ad \fs} \ra \subeq \Aut(\g)$ for the subgroup 
  generated by $e^{\ad \fs}$. We call $\fs$ {\it compactly embedded} if
  the group $\Inn_\g(\fs)$ has compact closure.
\item For $x \in \g$, we write $\cO_x := \Inn(\g)x$ for its adjoint orbit.   
\item If $\g$ is a Lie algebra, we write $\cE(\g)$ for the set of 
{\it Euler elements} $h \in \g$, i.e., $\ad h$ is non-zero and diagonalizable 
with $\Spec(\ad h) \subeq \{-1,0,1\}$. We call $h$ {\it symmetric} if 
$-h \in\cO_h$. The involution of $\g$ specified by~$h$
is denoted $\tau_h^\g := e^{\pi i \ad h}$. 
\item If $G$ is a group with Lie algebra $\g$ and $h \in \cE(\g)$, we
write $\tau_h \in \Aut(G)$ for the corresponding involution on $G$,
provided it exists, and $G_{\tau_h} := G \rtimes \{\id_G, \tau_h\}.$
\item Abstract Euler wedges:
  $\cG_E(G_{\tau_h}) = \{ (x,\sigma) \in
  \cE(\g) \times G\cdot \tau_h \:
  \sigma^2 =e, \Ad(\sigma) = \tau_x^\g\}.$
\end{itemize}

\medskip

{\bf Acknowledgments}: The work of V.M. is partly supported by
INdAM-GNAMPA, University of Rome Tor Vergata funding OANGQS CUP E83C25000580005, and the MIUR Excellence Department Project MatMod@TOV awarded to the Department of Mathematics, University of Rome Tor Vergata, CUP E83C23000330006.

\section{Orthogonal pairs of   Euler elements}
\mlabel{sec:2}

In this section we first extend some
of the results of \cite{MN21} for Euler elements in simple
Lie algebras to general Lie algebras, in particular
about the relation between the symmetry of $h$ and the
existence of orthogonal pairs $(h,k)$
    (Theorem~\ref{thm:sl2-gen}).
Then we describe a classification of all orthogonal pairs of Euler elements
  in simple Lie algebras under the action of the group
  $\Inn(\g)$ of inner automorphisms (Theorem~\ref{thm:1.4}). 
It is based on our classification of Euler elements 
in \cite{MN21} and on structural results on $3$-graded
simple Lie algebras and the fact that
any such pair generates a Lie subalgebra isomorphic
  to $\fsl_2(\R)$.

\subsection{Orthogonal pairs of Euler elements}

\begin{defn} \mlabel{def:diag-act}
  (a) We call a pair $(h,k)$ of Euler elements in a Lie algebra
  $\g$ {\it orthogonal} if $\tau_h^\g(k) = -k$.
We write $(h,k) \sim (h',k')$ if there exists a $\phi \in \Inn(\g)$ with
$(\phi(h),\phi(k)) = (h',k')$. 

\nin (b)   An Euler element $h \in \g$ is called {\it symmetric} if
  $-h \in \cO_h = \Inn(\g)h$. 
\end{defn}

\begin{lem} If an Euler element $h \in \g$ is symmetric,
    then $h \in [\g_1, \g_{-1}]$.   
\end{lem}

\begin{prf} First we note that
  $\fn := \g_1 + [\g_1,\g_{-1}] + \g_{-1} \trile \g$ is an ideal
  because it is invariant under $\ad \g_j$ for $j = 1,0,-1$.
  In the quotient Lie algebra $\fq := \g/\fn$, the image $\oline h$ of
  $h$ is central. On the other hand, there exists an inner automorphism
  $\phi \in \Inn(\g)$ with $\phi(h) = -h$, and this property is inherited
  by $\oline h$. This implies that $\oline h =  0$, i.e., $h\in
  \fn \cap \g_0 = [\g_1,\g_{-1}]$.   
\end{prf}

\begin{ex}
  \mlabel{ex:sl2-a} {\rm(Euler elements in $\fsl_2(\R)$)} 
In $\g= \fsl_2(\R)$ we have the orthogonal Euler elements 
\begin{equation}
  \label{eq:h0k0}
 h_0 = \frac{1}{2}\pmat{1 & 0 \\ 0 &-1} \quad \mbox{ and } \quad 
 k_0 := \frac{1}{2}\pmat{0 & 1 \\ 1 &0}.
\end{equation}
The group $\Inn(\g)$ acts transitively on the set 
  \[\cE(\fsl_2(\R))
    = \big\{ x \in \fsl_2(\R) \: \det(x) = -{\textstyle\frac{1}{4}}\big\} = 
    \Big\{ \pmat{a & b \\ c & -a} \: a^2 + bc = \frac{1}{4}\Big\},\]
which is diffeomorphic to the $2$-dimensional de Sitter space
$\dS^2$ (a one-sheeted hyperboloid).
From
  \[ \tau_{h_0}\pmat{a & b \\ c & -a}
    = \pmat{a & -b \\ -c & -a} \]
  it follows that the Euler elements orthogonal to $h_0$ are those of the
  form $k = \frac{1}{2} \pmat{0 & b \\ b^{-1} & 0},$ 
  and we obtain $k_0$ for $b = 1$. So, for $\fsl_2(\R)$,
    orthogonality of Euler elements is equivalent to their orthogonality
    with respect to the Cartan--Killing form.
  In view of \cite[Lemma~3.6]{MN21}, we also have
  $\tau_{k_0}(h_0) = -h_0$. 

Now it is easy to see that the pairs $(h_0, k_0)$ and $(h_0, -k_0)\sim (k_0,h_0)$ represent the two conjugacy classes of orthogonal pairs in $\fsl_2(\R)$.
For 
\begin{equation}
  \label{eq:z0}
  z_0 :=\frac{1}{2} \pmat{0 & 1 \\ -1 & 0} =  [h_0, k_0]
  \in   C_0 
\end{equation}
we then have
    \begin{equation}
      \label{eq:z0-rot}
e^{-\frac{\pi}{2} \ad z_0} h_0 = k_0 \quad \mbox{ and  } \quad        
e^{-\pi \ad z_0} h_0 =  e^{\pi \ad z_0} h_0 =  -h_0.
    \end{equation}
The subset 
  \begin{equation}
    \label{eq:c0insl2}
C_0 :=  C_{h_0, k_0} :=
    \Big\{ \pmat{a & b \\ c & -a} \: a^2 + bc \leq 0,
    b \geq 0 \Big\}
  \end{equation}
  is a pointed, generating, closed convex $\Inn(\g)$-invariant cone
  containing $z_0 = [h_0, k_0]$.
\end{ex}

The following proposition generalizes some results  from
\cite{MN21} on orthogonal Euler pairs. 

\begin{prop} \mlabel{prop:MN21-3-13} 
The following assertions hold:
  \begin{itemize}
  \item[\rm(a)] If $\g$ is semisimple and $(h,k)$ is an orthogonal
    pair of Euler elements, then $k$ is symmetric.
  \item[\rm(b)] If $\g$ is simple, then a
    pair $(h,k)$ of Euler elements in $\g$ is
  orthogonal if and only if $(k,h)$ is orthogonal. 
\item[\rm(c)] If $\g$ is semisimple and
  $(h,k)$ and $(k,h)$ are orthogonal pairs of Euler elements, then
  the Lie subalgebra generated by $h$ and $k$ is
  isomorphic to $\fsl_2(\R)$.
  \end{itemize}
\end{prop}

\begin{prf} (a)  We write $\g$ as a sum $\g = \g_1 \oplus \cdots \oplus \g_N$ 
of simple ideals~$\g_j$
 and $h = \sum_j h_j$, $k = \sum_j k_j$ with $h_j, k_j \in \g_j$.
 Then either $h_j = 0$, or $h_j$ is an Euler element in $\g_j$.
 The relation $\tau^\g_h(k) = -k$
  implies $\tau^\g_{h_j}(k_j) = -k_j$ for each $j$.
  If $k_j \not=0$, this implies $h_j \not=0$ and,
  by the simple case (\cite[Thm.~3.13]{MN21}), that $k_j$ is symmetric,
  hence that $k = \sum_j k_j$ is symmetric.

  \nin (b) follows from \cite[Thm.~3.13]{MN21}.
\begin{footnote}{The proof of \cite[Thm.~3.13]{MN21}
    is somewhat short in the complex case, where it is asserted that
    $h$ is contained in a certain subalgebra $\fs \cong \fsl_2(\C)^r$
    specified by a maximal system of strongly orthogonal roots.
    This issue is discussed in detail in \cite[Cor.~3.5]{MNO23}.   
}\end{footnote}

  \nin (c) We write $\g$, $h$ and $k$ as in (a).
  Then our assumption implies that $h_j \not=0$ is equivalent to
  $k_j \not=0$. We may thus assume w.l.o.g.\ that all $h_j$
  are non-zero.
    \cite[Thm.~3.13]{MN21} now implies the existence of Lie algebra
  homomorphisms
  \[ \phi_j \: \fsl_2(\R) \to \fs_j, \quad \mbox{ with } \quad
    \phi_j(h_0) = h_j, \quad \phi_j(k_0) = k_j,\] 
  with the notation from Example~\ref{ex:sl2-a}. 
  Then $\phi := \sum_j \phi_j \: \fsl_2(\R) \to \fs \subeq \g$
  also is a homomorphism of Lie algebras
  satisfying $\phi(h_0) = h$ and $\phi(k_0) = k$. 
\end{prf}

\begin{rem} In Proposition~\ref{prop:MN21-3-13}(b),
    the simplicity assumption is crucial.
    If, e.g., $\g = \g_1 \oplus \g_2$ is a direct sum
    and $h = h_1 + h_2$ is an Euler element with $h_1, h_2 \not=0$
    and $k \in \g_2$ is an Euler element satisfying
    $\tau_h(k) = - k$, then $(h,k)$ is orthogonal but
    $(k,h)$ is not. 
\end{rem}

\begin{thm} {\rm(Symmetry and orthogonal pairs)} \mlabel{thm:sl2-gen}
For an Euler element  $h \in \g$, the following   are equivalent:
  \begin{itemize}
  \item[\rm(a)] $h$ is symmetric.
  \item[\rm(b)] There exists an Euler element $k \in \cE(\g)$
        such that $(h,k)$ and $(k,h)$ are orthogonal.
 \end{itemize}
 For $k \in \cE(\g)$, condition {\rm(b)} is equivalent to:
 \begin{itemize}
 \item[\rm(c)] The Lie algebra
   generated by $h$ and $k$ is isomorphic to $\fsl_2(\R)$.
 \end{itemize}
      
\end{thm}

\begin{prf} (a) $\Rarrow$ (b):
  Suppose first that $h$ is symmetric.
  In view of \cite[Prop.~3.2]{MN21}, there exists a
  Levi decomposition $\g = \fr \rtimes \fs$ with $h \in \fs$
  and $h$ is a symmetric Euler element in $\fs$.
  Write $\fs = \fs_1 \oplus \cdots \oplus \fs_N$ with simple
  ideals $\fs_j$ and $h = \sum_j h_j$ with $h_j \in \fs_j$.
  Then either $h_j = 0$, or $h_j$ is an Euler element in $\fs_j$.
  In the latter case \cite[Thm.~3.13]{MN21} implies the existence
  of an Euler element $k_j \in \fs_j$ such that
  $(h_j, k_j)$ and $(k_j, h_j)$ are orthogonal.
  Then $k := \sum_j k_j \in \fs$ is an Euler element in $\fs$
  for which   $(h,k)$ and $(k,h)$ are orthogonal.
  Since $h$ and $k$ generate a subalgebra isomorphic to $\fsl_2(\R)$
  (Proposition~\ref{prop:MN21-3-13}(c)),
  $h$ and $k$ are conjugate by Example~\ref{ex:sl2-a}.
  Hence $k$ is also an Euler element of $\g$. 
\nin (b) $\Rarrow$ (a):  
  Let $h,k \in \cE(\g)$ be such that both $(h,k)$ and $(k,h)$ are orthogonal.
We write  $q_\fs \: \g \to \fs := \g/\fr$ for the projection
  onto the quotient by the solvable radical $\fr \trile \g$.
  In particular $\fs$ is isomorphic to any Levi subalgebra of $\g$.

  \nin {\bf Step 1:} We claim that $h,k \not\in\fr$.
  We argue by contradiction.   If $k \in \fr$, then
  $\tau^\g_h(k) = -k$ implies that $k \in \fr_1(h) + \fr_{-1}(h)
  \subeq [\g,\fr]$, so that $\ad k$ is nilpotent; a contradiction
  to the Euler element  property. As our assumptions on $h$ and $k$
  are symmetric, we also have $h\not\in \fr$.
Thus $h_\fs := q_\fs(h)$ and $k_\fs := q_\fs(k)$ are Euler elements in
$\fs$ for which $(h_\fs, k_\fs)$ and $(k_\fs, h_\fs)$ are orthogonal.
Both are symmetric by Proposition~\ref{prop:MN21-3-13}(a),(b), and
Proposition~\ref{prop:MN21-3-13}(c) implies that the
subalgebra they generate is isomorphic to $\fsl_2(\R)$. 

\nin {\bf Step 2:} We may w.l.o.g.\ assume that
the Lie algebra $\g$ is generated by $h$ and $k$.
Then Step~1 implies that $\fs \cong \fsl_2(\R)$.
Moreover, the orthogonality of $(h,k)$ and $(k,h)$ implies that
\[ k = [h,[h,k]] \quad \mbox{ and } \quad h = [k,[k,h]].\]
Therefore $\g = [\g,\g]$. This in turn implies that
 $\fr = [\g,\fr]$ is a nilpotent ideal.

 \nin {\bf Step 3:}
 According to \cite[Prop.~I.2]{KN96}, applied to the subgroup $\exp(\R \ad k)
 \subeq \Aut(\g)$,
there exists a Levi subalgebra $\fs \subeq \g$ with
$[k,\fs] \subeq \fs$. We then have
\[ \g = \fr \rtimes \fs, \quad
  k = k_\fr + k_\fs, \quad  h = h_\fr + h_\fs.\]
The relation $[k,\fs] \subeq \fs$ and $[k,k] = 0$ imply that
$[k,k_\fr] = [k,k_\fs] = 0$. Thus
$k_\fr = k - k_\fs$ is a difference of two $\ad$-diagonalizable
commuting elements, hence diagonalizable. It is also contained
in the nilpotent ideal $\fr$, hence $\ad$-nilpotent, and therefore
central. So $k_\fs$ is also Euler in $\g$ and
$\tau^\g_k = \tau^\g_{k_\fs}$. We conclude that
\begin{equation}
  \label{eq:dag1}
 \tau_k^\g(h_\fs) = \tau^\g_{k_\fs}(h_\fs)
 = q(\tau^\g_k(h)) = q(-h) = - h_\fs.
\end{equation}
As $h_\fs \in \Inn_\g(\fs) k_\fs$ holds in $\fs \cong \fsl_2(\R)$
(Example~\ref{ex:sl2-a}), $h_\fs$ also is a symmetric Euler element in~$\g$. 

\nin {\bf Step 4:} We consider the solvable Lie subalgebra
\[ \fb := \fr + \R h = \fr + \R h_\fs \cong \fr \rtimes_{\ad h} \R.\]
Our strategy is to show that $h = h_\fs$, and since $h_\fs \in \fs$
is a symmetric Euler element in $\fs$, this implies (a). We claim that
\begin{equation}
  \label{eq:ddag}
  \cE(\g) \cap \fb = \Inn_\g(\fb)(\pm h + \fz(\g)).
\end{equation}
The inclusion ``$\supeq$'' is trivial. 
To prove the converse, we first observe that
\[ \fh := \fb_0(h)  = \fr_0(h) + \R h \subeq \fb \]
is a Cartan subalgebra of $\fb$. In fact, as $h$ is central in $\fh$ 
and $\fr$ is nilpotent, $\fh$  is nilpotent, and the
eigenspace decomposition of $\ad h$ on $\fb$ implies that $\fh$ is
self-normalizing (\cite[Prop.~6.1.6]{HN12}). Now let $m \in \fb \cap \cE(\g)$.
Them $m$ is $\ad$-semisimple, hence contained in a
Cartan subalgebra (\cite[Prop.~6.1.12]{HN12}). As all Cartan subalgebras
of the solvable Lie algebra $\fb$ are conjugate under
inner automorphisms (\cite[Ch.~7, \S 3, no.~4, Thm.~3]{Bo90b}), this implies
that $\Inn(\fb)m$ intersects $\fh$. We may therefore assume that
$m \in \fh \cap \cE(\g)$. We write
$m = \lambda h + z$ with $z \in \fr_0(h)$. Then
$\ad z$ is nilpotent, and since the eigenvalues of $\ad m$ are
$\pm 1$, we obtain $\lambda \in \{ \pm 1\}$.
Further $z = m - \lambda h$ is $\ad$-nilpotent and $\ad$-diagonalizable,
hence central in $\g$. This shows that
$m \in \pm h + \fz(\g)$, and \eqref{eq:ddag} follows. Note that
$\tau_h = \tau_m$. 

\nin {\bf Step 5:} We consider the adjoint orbit
$\cO^B_m := \Inn(\fb)m = \Inn(\fr)m$.

Let $G \cong R \rtimes S$ be a simply connected Lie group
with Lie algebra $\g$. As $\fr$ is nilpotent
and $R$ simply connected, the exponential function
$\exp_R \: \fr \to R$ is a diffeomorphism (\cite[Thm.~11.2.10]{HN12}).
This implies easily that
\[ R^m = R^h = R^{\tau_h} = R^{\tau_m} = \exp(\fr_0(h))
  \quad \mbox{ and } \quad
  R^{-\tau_h} := \{ g \in R \: \tau_h(g) = g^{-1} \}
  = \exp(\fr^{-\tau^\g_h}).\] 
Since the map 
\[ \cO_m^B \cong R/R^{\tau_h} \to R^{-\tau_h}, \quad
\Ad(g)m \mapsto   g R^{\tau_h} \mapsto g \tau_h(g)^{-1} \]
is a diffeomorphism, mapping
$e^{\ad x}m$ to $\exp(2x)$ for $x \in \fr^{-\tau_h^\g}$, it follows that
\[ \Phi \: \fr^{-\tau_h^\g} = \fr_1(h) + \fr_{-1}(h)
  \to \cO_m^B, \quad x \mapsto e^{\ad x}h \]
is a diffeomorphism.
Replacing $\fr$ by a quotient Lie algebra
$\oline \fr := \fr/\R z$, where $z \in \fz(\g)$, 
we likewise obtain a diffeomorphism
\begin{equation}
  \label{eq:dddag}
 \oline\Phi \: \oline\fr^{-\tau_h^\g} \cong \fr_1(h) + \fr_{-1}(h)
 \to \Inn(\oline\fr)m, \quad x \mapsto e^{\ad x}m.
\end{equation}

We now analyze the action of the involution $\tau_k^\g$ in this context.
By assumption $\tau^\g_k(h) = - h$, and we have also seen
in \eqref{eq:dag1} that $\tau^\g_k(h_\fs) = \tau^\g_{k_\fs}(h_\fs) = - h_\fs$. 
The first relation implies that the involutions
$\tau^\g_k$ and $\tau^\g_h$ commute, so that $\tau^\g_k$ leaves
$\fr^{-\tau^\g_h}$ invariant.

We pick  $m = h + z  \in h + \fz(\g)$ such that
$h_\fs \in \Inn(\fb)m$ (cf.\ Step 4). As we have seen above,
there exists a unique $x \in \fr^{-\tau^\g_h}$ with
$h_\fs = e^{\ad x}(h + z) = e^{\ad x}h + z$. Then
\begin{equation}
  \label{eq:dag4}
  - e^{\ad x}h - z = - h_\fs\ {\buildrel\eqref{eq:dag1}\over =}\ \tau^\g_k(h_\fs)
  = e^{\ad \tau^\g_k(x)} \tau^\g_k(h) + z
  = -e^{\ad \tau^\g_k(x)} h + z.
\end{equation}
Modulo $\R z$, this leads with \eqref{eq:dddag} to the relation
$e^{\ad x}h = e^{\ad \tau^\g_k(x)} h,$ 
and hence to $x = \tau^\g_k(x)$ because $\oline\Phi$ is a diffeomorphism.
Thus relation \eqref{eq:dag4} simplifies to
$- e^{\ad x}h - z   = -e^{\ad x}h + z,$ which leads to $z = 0$.
Therefore $h_\fs \in \Inn(\g)h$, and the fact that $h_\fs$ is a symmetric
Euler element in $\g$ (Step 3) shows that $h$ is symmetric.

\nin (b) $\Rarrow$ (c): We continue the argument with the
context from the preceding paragraph. 
The relations $h_\fs = e^{\ad x}h$ (Step 5) and $x \in \fr_0(k)$
further show that
\[ e^{\ad x}(h,k) = (h_\fs, k).\]
We may thus assume that $h = h_\fs$. Writing
$k = k_\fs + k_\fr$ with $k_\fr \in \fz(\g)$ (Step 4), 
we then obtain
\[ -k_\fs - k_\fr = -k =
  \tau^\g_h(k) = 
  \tau^\g_{h_\fs}(k) = - k_\fs + k_\fr,\]
which leads to $k_\fr = 0$. This shows that $k = k_\fs$,
and since $\g$ is generated by $h$ and $k$, we eventually obtain $\g = \fs$.

\nin (c) $\Rarrow$ (a) follows from the fact that all Euler elements in
$\fsl_2(\R)$ are symmetric (Example~\ref{ex:sl2-a}). 
\end{prf}

\subsection{Euler elements in simple Lie algebras} 
  
The classification of the $\Inn(\g)$-orbits  of symmetric Euler elements is
reproduced below from \cite[Thm.~3.10]{MN21}
(see also \cite{KA88, Kan98, Kan00}).

\begin{theorem} \mlabel{thm:classif-symeuler}
{\rm(\cite[Thm.~3.10]{MN21})} Suppose that $\g$ is a non-compact simple 
real Lie algebra, with restricted root system $\Sigma(\g,\fa)$. 
We follow the conventions of the tables in {\rm\cite{Bo90a}}
for the classification of irreducible root systems and the enumeration 
of the simple roots $\Pi = \{\alpha_1, \ldots, \alpha_n\}$.
We write $h_1, \ldots, h_n$ for the dual basis of $\fa$, i.e.,
$\alpha_k(h_j) = \delta_{kj}$. Then every Euler element $h \in \fa$ on which 
$\Pi$ is non-negative is one of  $h_1, \ldots, h_n$, and, for 
every irreducible root system, the Euler elements among the $h_j$ are 
 the following, which represent the $\Inn(\g)$-orbits: 
\begin{align} 
&A_n: h_1, \ldots, h_n, \quad 
\ \ B_n: h_1, \quad 
\ \ C_n: h_n, \quad \ \ \ D_n: h_1, h_{n-1}, h_n, \quad 
E_6: h_1, h_6, \quad 
E_7: h_7.\label{eq:eulelts2}
\end{align}
For the root systems $BC_n$, $E_8$, $F_4$ and $G_2$ no Euler element exists.\begin{footnote}{If $\g$ is a real or complex
    split Lie algebra, then all automorphisms of the Dynkin diagram
    extend to automorphisms of $\g$.
    Therefore the Euler elements $h_j$ and $h_{n-j+1}$ are conjugate under
    $\Aut(\g)$ for $A_n$; likewise are $h_{n-1}$ and $h_n$ for $D_n$,
    and $h_1$ and $h_6$ for $E_6$.  
  }\end{footnote}

The symmetric Euler elements are 
\begin{equation}
  \label{eq:symmeuler}
A_{2n-1}: h_n, \qquad 
B_n: h_1, \qquad C_n: h_n, \qquad 
D_{2n}: h_1, h_{2n-1},h_{2n}, \quad 
D_{2n+1}: h_1, \qquad 
E_7: h_7.  
\end{equation}
For type $D_{2n}, n \geq 2$, there are three conjugacy classes of symmetric
Euler elements,  
and in all other cases there is only one.
\end{theorem}

\begin{rem}
If $\g$ is hermitian, then 
$\Sigma(\g,\fa)$ is of type $C_n$ or $BC_n$ and
Euler elements only exist in the first case, 
which means that $\g$ is of tube type. In this case 
$\Inn(\g)$ acts transitively on $\cE(\g)$
(\cite[Prop.~3.11]{MN21}). 
\end{rem}

For Table 1 below, we recall the Lie algebras
\[ \so^*(2n) := \fu_n(\bH,i\1) :=\{ x \in \gl_n(\H) \:
  x^* i  + i x = 0 \}.\]
As in Example~\ref{ex:complexcase} below, we write $\g$
  for a complex simple Lie algebra and $\g^\circ$ for a
  hermitian real form of $\g$.  \\

\begin{tabular}{||l|l|l|l|l||}\hline
{} $\g^\circ$ \mbox{(hermitian)}  & $\Sigma(\g^\circ, \fa^\circ)$ & $\g = (\g^\circ)_\C$ & $\Sigma(\g,\fa)$ & {\rm symm.\ Euler element}\  $h\in \g^\circ$ \\ 
\hline\hline 
$\su_{n,n}(\C)$ & $C_n$ & $\fsl_{2n}(\C)$ & $A_{2n-1}$ & $h_n$ \\ 
 $\so_{2,2n-1}(\R), n > 1$ & $C_2$ & $\so_{2n+1}(\C)$ & $B_n$ & $h_1$ \\ 
$\sp_{2n}(\R)$ & $C_n$ & $\fsp_{2n}(\C)$ & $C_n$ & $h_n$ \\ 
 $\so_{2,2n-2}(\R), n > 2$ & $C_2$ & $\so_{2n}(\C)$ & $D_{n}$ & $h_1$ \\ 
  $\so^*(4n) = \fu_{2n}(\bH,i\1) $ & $C_n$ & $\so_{4n}(\C)$  & $D_{2n}$ & $h_{2n-1},  h_{2n}$ \\ 
$\fe_{7(-25)}$ & $C_3$ & $\fe_7$ & $E_7$ & $h_7$ \\ 
\hline
\end{tabular} \\[2mm] {\rm Table 1: Simple hermitian Lie algebras $\g^\circ$ 
  with Euler element, i.e., of tube type}\\

\subsection{The classification of orthogonal Euler pairs} 
\mlabel{subsec:2.2}

To classify orthogonal pairs of Euler elements, we have to determine, 
for a given symmetric Euler element 
$h\in \cE(\g)$ (cf.\ Proposition~\ref{prop:MN21-3-13}(a),(b))
and $G = \Inn(\g)$,
the orbits of the stabilizer group $G^h$ of $h$ in
the $G^h$-invariant subset $\cE(\g) \cap \fq$, where 
\begin{equation}
  \label{eq:fq}
  \fq := \g^{-\tau_h^\g} = [h,\g] = \g_1(h) +\g_{-1}(h).
\end{equation}
Below we also write 
\begin{equation}
  \label{eq:fh}
 \fh := \g^{\tau_h^\g} = \g_0(h), \quad \mbox{ so that } \quad 
 \g = \fh \oplus \fq.
\end{equation}

\begin{thm} \mlabel{thm:1.4} {\rm(Classification Theorem for pairs 
of orthogonal Euler elements)} 
Let $\g$ be a simple real Lie algebra and 
$h \in \cE(\g)$ a symmetric Euler element. 
Let $\theta$ be a Cartan involution with $\theta(h) = -h$ and 
$\fa \subeq \fp := \g^{-\theta}$ be a maximal abelian subspace  
containing~$h$, so that $\fa \subeq \fh_\fp$. In 
\[ \Sigma_1 := \{ \alpha \in \Sigma(\g,\fa) \: \alpha(h) = 1\} \]
we pick a  maximal set $\{ \gamma_1, \ldots, \gamma_r\}$ of long strongly 
orthogonal roots.
\begin{footnote}{Two roots $\alpha,\beta$ are called {\it strongly orthogonal}
if neither $\alpha + \beta$ nor $\alpha - \beta$ is a root.
We refer to \cite[\S 3.2]{MNO23} for more on systems of strongly orthogonal
roots.}\end{footnote}
Then there exist elements $e_j \in \g_{\gamma_j}$ such that 
the Lie subalgebras 
\[ \fs_j := \Spann_\R \{e_j, \theta(e_j), [e_j, \theta(e_j)]\}, \quad 
j =1,\ldots, r, \] 
are isomorphic to $\fsl_2(\R)$.  
Then 
\begin{equation}
  \label{eq:defc}
  \fc := \Spann \{ k_1, \ldots, k_r \}
\quad \mbox{ for }  \quad k_j := e_j - \theta(e_j) 
\end{equation}
is maximal abelian in $\fq_\fp := \fq \cap \fp$ and contained in the 
subalgebra $\g^* := \fh_\fk \oplus \fq_\fp = \Fix(\theta \tau_h^\g)$,
so that we have an inclusion 
of restricted root systems $\Sigma(\g^*,\fc) \subeq \Sigma(\g,\fc)$. 
We consider the Euler elements 
\begin{equation}
  \label{eq:hj}
 k^j := k_1 + \cdots + k_j - k_{j+1} - \cdots - k_r \in \fc,  
 \quad j = 0,\ldots, r. 
\end{equation}
Then $\Sigma(\g,\fc)$ is of type $C_r$ and 
$\Sigma(\g^*,\fc)$ is of type $A_{r-1}$, $C_r$ or $D_r$,
as listed in {\rm Table 2} below.
All $G_e^h$-orbits in $\cE(\g) \cap \fq$ are $G^h$-invariant, 
 and a set of representatives is given by: 
\begin{itemize}
\item[\rm(A)] $k^0, \ldots, k^r$, if $\Sigma(\g^*,\fc)$ is of type $A_{r-1}$ and
  $r \geq 1$. \begin{footnote}
    {The case $r = 1$ occurs only for $\g = \su_{1,1}(\C) \cong \fsl_2(\R)$
      in Table 2 below. Then $\g^* \cong \R$ is of type $A_0$.}
  \end{footnote}
\item[\rm(C)] $k^r$, if $\Sigma(\g^*,\fc)$ is of type $C_r, r \geq 2$. 
\item[\rm(D)] $k^{r-1}, k^r$, if $\Sigma(\g^*,\fc)$ is of type $D_r, r \geq 2$.
\end{itemize}
Accordingly, the pairs $(h,k^j)$ represent the $\Inn(\g)$-conjugacy classes of
pairs of orthogonal Euler elements.
\end{thm}

Since $h$ is a symmetric Euler element and the root system 
$\Sigma(\g,\fa)$ is irreducible, $h$ corresponds to some 
$h_j$ in the list \eqref{eq:symmeuler} in Theorem~\ref{thm:classif-symeuler}. 
Geometrically, the symmetric Lie algebra $(\g^*,\tau_h^\g\res_{\g^*})$
  corresponds to a maximal negatively curved Riemannian subspace
  of the symmetric space $G/G^h$, 
  corresponding to the symmetric Lie algebra $(\g,\tau_h^\g)$.
  The tangent space of $G/G^h$ in the base point is $\fq$,
  and the subspace corresponding to $G^*/G^{*,h}$ is
  $\fq_\fp$. Its elements generate non-compact geodesics, whereas
  the elements of $\fq_\fk$ generate geodesics with compact closure.  
  The pair $(\g,\g^*)$ is non-compactly causal in the terminology of
  \cite{HO97, MNO23}.

\begin{ex} \mlabel{ex:sl2} {\rm(Euler pairs in $\fsl_2(\R)$)}
  We continue the discussion from Example~\ref{ex:sl2-a}
  of Euler elements in $\g = \fsl_2(\R)$, where
  $(h_0, k_0)$ and $(h_0, -k_0)\sim (k_0,h_0)$
represent the two conjugacy classes of orthogonal pairs. 
The Lie subalgebra $\g^* = \fh_\fk + \fq_\fp  = \fq_\fp = \fc\cong \R$
is abelian, so that $\Sigma(\g^*,\fc) = \eset$ (type $A_0$) and the two Euler elements 
$k^0$ and $k^1 = - k^0$ in $\fc$ are not conjugate under $G^h_e \cong \R$, 
acting by positive dilations. 

For two orthogonal Euler elements $h,k$,
  let $C_{h,k} \subeq \fsl_2(\R)$ be the pointed, generating 
closed convex invariant cone generated by $[h, k]$ (cf.~\eqref{eq:c0insl2}). 
Then $C_{h, k} = - C_{k,h}$, so that
the two conjugacy classes of orthogonal 
  pairs of Euler elements $(h,k)$ are distinguished by 
the corresponding cone.
\end{ex} 

\begin{ex} {\rm(The complex case)} \mlabel{ex:complexcase}
  If $\g$ is complex simple, containing a pair $(h,k)$ of orthogonal
  Euler elements (Example~\ref{ex:sl2}),
  then $h$ is symmetric (Proposition~\ref{prop:MN21-3-13}).
  This implies that $\fg \cong (\g^\circ)_\C$, 
  where $\g^\circ$ is hermitian of tube type
  (cf.\ \cite[Prop.~3.11]{MN21}).  
Let $\g^\circ = \fk^\circ \oplus \fp^\circ$ be a Cartan decomposition,
so that the associated Cartan decomposition of $\g$ is 
$\fk \oplus \fp$ with 
$\fk = \fk^\circ \oplus i \fp^\circ$ and $\fp = i \fk$. 
If $\ft^\circ \subeq \fk^\circ$ is a Cartan subalgebra, then 
$\fa := i \ft^\circ \subeq \fp$ is maximal abelian and the restricted 
root system $\Sigma(\g,\fa)$ is  of type 
$A_{2n-1}, B_n, C_n, D_n$ or $E_7$ 
(cf.~Table 1). 
Only for $D_{2n}$, corresponding to 
$\g^\circ = \so^*(4n)$, there are three conjugacy classes of symmetric 
Euler elements, but in all other cases they are unique 
(Theorem~\ref{thm:classif-symeuler}). 
In the $D_{2n}$-case, $h_{2n-1}$ and $h_{2n}$ are conjugate under a diagram automorphism. 

The Euler element $h \in \fa = i \ft^\circ \subeq i \fk^\circ$ 
is contained in the one-dimensional space $i \fz(\fk^\circ)$. 
This shows that $\fh = \g_{\C,0}(h)= \fk^\circ_\C$,
and $\fq = \fp^\circ_\C$, which in turn implies 
that $\g^* = \g^\circ$. For the 
maximal abelian subspace $\fc \subeq \fp^\circ = \fq_\fp$, 
the restricted root system $\Sigma(\g^\circ, \fc)$ is of 
type~$C_r$ (cf.~\cite[Thm.~XII.1.14]{Ne99}). 

As $\g \cong\g^\circ_\C$, the root systems 
$\Sigma(\g,\fc)$ and $\Sigma(\g^\circ,\fc)$ coincide, hence are both of type 
$C_r$. We conclude with Theorem~\ref{thm:1.4} that 
there is only one $G^h_e$-orbit in $\fq \cap \cE(\g)$. 
It corresponds to the unique class of Euler elements 
associated to the root system $C_r$ 
(Theorem~\ref{thm:classif-symeuler}). 

This correspondence shows that conjugacy classes of 
orthogonal pairs $(h,k)$ of Euler elements in the complex simple
Lie algebra $\g$
correspond to $G$-orbits in $\cE(\g)$, resp., to hermitian tube type real forms
of $\g$. For $\Sigma(\g,\fa)$ not of type $D_{2n}$, 
there is a single orbit (Theorem~\ref{thm:classif-symeuler}),
and, for type $D_{2n}$, there are three orbits, 
one corresponding to the hermitian real form
$\g^\circ \cong \so_{2,4n-2}(\R)$, and two 
corresponding to the hermitian real form
$\g^\circ \cong \so^*(4n)$. The latter two 
are conjugate under an outer automorphism. 
\end{ex}

The following table has been composed with the aid of Table 5 in 
\cite[p.~135]{Kan00}. It identifies in all cases the data mentioned 
in the theorem. For two cases the $h$-column contains two
  entries $h_{2r-1}$ and $h_{2r}$. Then the Lie algebra $\g$ is split
  (all root spaces are $1$-dimensional) with
  root system $D_{2r}$, so that the diagram reflection extends to a
  Lie algebra automorphism exchanging these two Euler elements. Therefore
  both lead to isomorphic configurations.
\\

 \hspace{-16mm}
\begin{tabular}{||l|l|l|l|l|l|l||}\hline
 $\g$   & $\g^* = \Fix(\tau_h^\g\theta)$\phantom{\Big(} & $\Sigma(\g,\fa)$  & $h$ & $\g_1(h)$ & 
$\Sigma(\g,\fc)$ &  $\Sigma(\g^*,\fc)$  \\ 
\hline\hline 
Complex type \phantom{\Big)} &&&& &&  \\
\hline 
 $\fsl_{2r}(\C)$ & $\su_{r,r}(\C)$ &  $A_{2r-1}$ & $h_r$ & $M_r(\C)$  & $C_r$ & $C_r$ \\
 $\sp_{2r}(\C)$ & $\sp_{2r}(\R)$ &  $C_{r}$ & $h_r$ & $\Sym_r(\C)$   & $C_r$  & $C_r$ \\
 $\so_{2n+1}(\C), n > 1$ & $\so_{2,2n-1}(\R)$ &   $ B_{n}$ & $h_1$ & $\C^{2n-1}$   & $C_2$ & $C_2$ \\
 $\so_{2n}(\C), n > 2$ & $\so_{2,2n-2}(\R)$ &  $ D_{n}$ & $h_1$ & $\C^{2n-2}$   & $C_2$  & $C_2$ \\
  $\so_{4r}(\C)$ & $\so^*(4r)$
 & $ D_{2r}$ & $h_{2r-1},h_{2r}$ & $\Skew_{2r}(\C)$   & $C_r$ & $C_r$ \\
$\fe_7(\C)$ & $\fe_{7(-25)}$ & $E_7$ & $h_7 $ & $\Herm_3(\bO)_\C$   & $C_3$ & $C_3$  \\
\hline 
Cayley type (CT) \phantom{\Big)}&& &&&&  \\
\hline 
  $\su_{r,r}(\C), r \geq 1$ & $\R \oplus \fsl_r(\C)$
      &$C_{r}$ & $h_r$ & $\Herm_r(\C)$  &  $C_r$ & $A_{r-1}$\\
$\sp_{2r}(\R), r > 1$ & $\R \oplus \fsl_r(\R)$ & $C_{r}$ & $h_r$ & $\Sym_r(\R)$   & $C_r$ &  $A_{r-1}$ \\
 $\so_{2,n+1}(\R), n > 1$ & $\R \oplus \so_{1,n}(\R)$ & $C_2$ & $h_2$ & $\R^{1,n}$   & $C_2$ & $A_1$ \\
 $\so^*(4r), r > 1$ &$\R \oplus \fsl_r(\bH)$ & $C_r$ & $h_r$ & $\Herm_r(\H)$   & $C_r$  & $A_{r-1}$ \\
 $\fe_{7(-25)}$ & $\R \oplus \fe_{6(-26)}$ &  $C_3$ & $h_3$ & $\Herm_3(\bO)$   & $C_3$  & $A_2$\\
\hline
Split type (ST) \phantom{\Big (} &&&&&&  \\
\hline 
 $\fsl_{2r}(\R)$, $r > 1$ & $\so_{r,r}(\R)$ & $A_{2r-1}$ & $h_r$ & $M_r(\R)$  & $C_r$ & 
$D_r$ \\
 $\so_{2r,2r}(\R), r > 1$ & $\so_{2r}(\C)$ & $D_{2r}$ & $h_{2r-1},h_{2r}$ & $\Skew_{2r}(\R)$   & $C_r$ & $D_r$ \\
 $\fe_7(\R)$ & $\fsl_4(\bH)$ &  $E_7$ & $h_7 $ & $\Herm_3(\bO_{\rm split})$   & $C_3$ & $D_3$ \\
  $\so_{p+1,q+1}(\R)$ & $\so_{1,p}(\R) \oplus$
                 & $ B_{p+1}\, (p<q)$ & $h_1$ & $\R^{p,q}$ & 
                                                             $C_2$ & $D_2 \cong$  \\
$p,q >1$ & \ \ \ \  \ \ $ \so_{1,q}(\R)$& $ D_{p+1}\, (p = q)$&&&& \ \ \ $A_1^{\oplus 2}$ \\ 
\hline 
Non-split type \phantom{\Big)} &&& && & \\[-3mm]
(NST) \phantom{\Big)} &&& && & \\
\hline 
 $\fsl_{2r}(\H)$ & $\fu_{r,r}(\bH)$ &  $A_{2r-1}$ & $h_r$ & $M_r(\H)$  & $C_r$ &  $C_r$ \\
 $\fu_{r,r}(\H)$ & $\sp_{2r}(\C)$ &  $C_{r}$    & $h_r$ & $\Aherm_{r}(\H)$ & $C_r$& $C_r$  \\
 $\so_{1,d+1}(\R), d \geq 2$ & $\so_{1,d}(\R)$ & $A_1$ & $h_1$ & $\R^{0,d}$ & $C_1$ & $C_1$ \\
\hline
\end{tabular} \\[2mm] {\rm Table 2: Simple $3$-graded Lie algebras with 
  {\bf symmetric} Euler elements as in \eqref{eq:symmeuler}.}\\

\begin{prf} (of Theorem~\ref{thm:classif-symeuler}) 
For the existence of the set of long strongly 
orthogonal roots we refer to \cite[p.13]{Ta79} or \cite[p.~134]{Kan00}. 
From these references we infer the existence of elements 
$e_j \in \g_{\gamma_j}$ such that the subalgebras 
$\fs_j$ are isomorphic to $\fsl_2(\R)$.
In general $\dim \g_{\gamma_j} > 1$, as the
  example $\g = \so_{1,d+1}(\R)$ with $d$-dimensional root spaces shows.
We normalize $e_j$ in such a way that  
$k_j = e_j - \theta(e_j)$ is an Euler element of $\fs_j$.
Then loc.\ cit.\ further implies that 
$\fc$, as in \eqref{eq:defc}, is maximal abelian in $\fq_\fp$. 
Hence, by \cite[Cor.~III.9]{KN96}, 
every Euler element $x \in \fq$ is conjugate under 
$G^h_e = \Inn_\g(\g^{\tau_h^\g}) = \Inn_\g(\fh)$ to an element of~$\fc$. 

The restricted root system $\Sigma(\g,\fc)$ is of type $C_r$ 
(\cite[Prop.~II.2.1]{Kan00}, \cite[p.~596]{Kan98}): 
\[ C_r = \{ \pm 2 \eps_j, \pm \eps_j \pm \eps_\ell\: 1 \leq j \not= \ell
\leq r\}, \] 
where 
$\eps_j(k_m) = \frac{\delta_{jm}}{2}$ for $j,m = 1,\ldots, r$. 
Therefore $k \in \fc$ is an Euler element of $\g$ if and only if it can be written 
in the form 
\begin{equation}
  \label{eq:x}
 k := \sum_{j = 1}^r \sigma_j k_j 
 \quad \mbox{ with } \quad \sigma_j \in \{\pm 1\}.
\end{equation}

\nin {\bf Step 1 (orbits of $G^h_e$).} 
Which of these elements $k$ are conjugate under the connected 
group  $G^h_e = \Inn_\g(\g^{\tau_h^\g}) = \Inn_\g(\fh)$ depends 
on the Weyl group of the root system $\Sigma(\g^*,\fc)$. 
The key result is Kaneyuki's generalization of 
Sylvester's Law of Inertia 
(\cite[Thm.~II.2.5]{Kan00}, \cite{Kan98}), asserting the classification 
stated under (A), (C) and (D). Note that the root system  
$\Sigma(\g^*,\fc)$ cannot be of type $B_r$ because it is contained 
in $\Sigma(\g,\fc)$ of type $C_r$. 

In case (A), the Weyl group is the symmetric group $S_r$, which can be used to 
move the positive signs to the left. In case (C), the Weyl group contains 
all sign changes and we have a single orbit, and, in case (D), we have only 
even sign changes and thus two orbits (cf.\ Remark~\ref{rem:b3}(b)). 
This proves the stated classification of the $G^h_e$-orbits. 

\nin \nin {\bf Step 2 (orbits of $G^h$).} It remains to show 
that all $G^h_e$-orbits are invariant under the possibly larger group~$G^h$. 
For case (C) there is nothing to show because there is only 
one $G^h_e$-orbit. We consider the cases (A) and (D) separately. 

\nin {\bf Case (A) (hermitian Lie algebras/Cayley type spaces).} As we see in Table 2, the 
root system $\Sigma(\g^*,\fc)$ is of type $A_{r-1}, r \geq 1$, 
\begin{footnote}{We put $A_0= \eset$, which occurs for $\g^* \cong \R$
    and $\g \cong \fsl_2(\R)$; cf.~Remark~\ref{rem:3.11}(a).}\end{footnote}
if and only if 
$\g$ is hermitian (cf.\ \cite[p.~134]{Kan00}). 
Then $E := \g_1(h)$ is a simple euclidean Jordan algebra 
and $\g_0(h) = \g^{\tau_h^\g}$ is the Lie algebra 
of the structure group of~$E$ (cf.\ Appendix~\ref{app:1} for more details).
Further, 
the elements $e_j \in \g_{\gamma_j}$ form a 
Jordan frame $(e_1, \ldots, e_r)$ in~$E$ (cf.\ \cite{FK94}). 
An inspection of the rank $2$-case, corresponding to the 
Lorentz group acting on Minkowski space $E = \R^{1,d-1}$, shows that 
the subgroup of $G^h_e \cong \Inn_\g(\g_0(h))$, preserving the
subspace~$\fc$, coincides with the permutations of the basis 
elements $k_1, \ldots, k_r$. We therefore obtain the representatives 
$k^j, j = 0,\ldots, r$,  
of the $G^h_e$-orbits of Euler elements in $\fq$. 
Considering the projections onto $\g_1(h) \cong E$, these elements 
correspond to the 
involutions 
\[ e^j := e_1 + \cdots + e_j - e_{j+1} - \cdots - e_r, 
\quad j = 0,\ldots, r \]  
in the Jordan algebra $E$. 

In \cite[Lemma~3.7]{Ne18} we have seen that $G^h$ can be identified with a 
subgroup of the automorphism group $\Aut(E_+)$ of the open symmetric cone $E_+
\subeq E$ of invertible squares. 
In view of \cite[p.~57]{FK94}, this group has a polar decomposition 
\[ \Aut(E_+) = \Aut(E,e) \cdot \exp(L(E)),\] 
where $L(x)y = xy$ is the Jordan multiplication, and 
$\Aut(E,e)$ is the group of automorphisms of the unital 
Jordan algebra $(E,e)$. Since $\Aut(E,e)$ preserves all 
$r+1$ connected components of the set $E^\times$ of invertible elements   
(which are specified in terms of the signs of spectral values), 
we see that $G^h$ and $G^h_e$ have the same orbits in $E^\times$. 
In particular, the Euler elements $k^0, \ldots, k^r$ also represent 
the $G^h$-orbits in $\fq \cap \cE(\g)$. 

\nin {\bf Case (D) (split type)} This concerns only the four cases: 
\[ \fsl_{2r}(\R), \quad \so_{2r,2r}(\R), \quad \fe_7(\R), \quad 
\mbox{ and }\quad \so_{p,q}(\R), p,q > 2.\] 

\nin (a) For $\g = \fsl_{2r}(\R)$, we have $E = M_r(\R)$ 
and $G \cong \PSL_{2r}(\R) \cong \SL_{2r}(\R)/\{\pm \1\}$. 
Here
\[ G^h \cong {\rm S}(\GL_r(\R) \times \GL_r(\R))/\{\pm\1\} \]
acts on $E$ by 
$(g_1, g_2)x = g_1 x g_2^{-1}$ for 
$g_1, g_2 \in \GL_n(\R)$ with $\det(g_1) \det(g_2) = 1$. 

The set $E^\times = \GL_r(\R)$ of invertible elements in $E$
has two connected components, 
specified by the sign of the determinant.
Both are invariant  under the action of $G^h$.
Since $k^j$ corresponds to an element of determinant $(-1)^{r-j}$,
the Euler elements
$k^{r-1}$ and $k^r$ are not conjugate under~$G^h$. 

\nin (b) For $\g = \so_{2r,2r}(\R)$, we have $E = \Skew_{2r}(\R)$ 
and $G \cong \SO_{2r,2r}(\R)_e/\{\pm \1\}$. In \cite{MNO24} it is
shown that the group $G^h$ is connected, i.e., $G^h = G^h_e$.

\nin (c) For $\g = \fe_7(\R)$ (the split real form) and $E = \Herm_3(\bO_{\rm split})$,  
we know already that 
$\Sigma(\g^*,\fc)$ is of type $D_2$, so that we have at most two orbits. 
As $r = 3$, these orbits are represented by  $\pm e$, where 
$e \in E$ is the Jordan unit. It therefore suffices to show that 
$-e \not\in \Ad(G^h)e$. In Proposition~\ref{prop:pos}, 
we show that every $g \in G^h$ satisfies 
$\Ad(g)^*\detE \in \R_+ \detE$. Therefore 
\[ \detE(\Ad(g)e) \in\R_+ \detE(e) = \R_+ 
\quad \mbox{ and } \quad \detE(-e) = (-1)^3 = -1 \] 
imply that $\Ad(g)e \not= -e$.

\nin (d) For $\g = \so_{p,q}(\R)$, $p,q > 2$, we consider the 
matrix realization 
\[ \so_{p,q}(\R) 
  = \bigg\{ \pmat{ a & b \\ b^\top & d} \: a \in \so_p(\R), d \in \so_q(\R),
  b \in M_{p,q}(\R) \bigg\}. \]
Then  we obtain for $1 \leq j \leq p < m \leq p+q$
  (boost) Euler elements $h_{j,m}$ by 
  \begin{equation}
    \label{eq:hjk}
h_{j,m}\be_j = \be_m, \quad h_{j,m} \be_m = \be_j \quad \mbox{ and } \quad 
h_{j,m}\be_\ell = 0 \quad \mbox{ for } \quad  \ell \not=j,m.
  \end{equation}
We consider the Euler element $h := h_{1,p+q}$.
In the Jordan algebra $E = \R^{p-1,q-1}$, the subset 
\[ E^\times = \{ (x,y) \in \R^{p-1} \times \R^{q-1} \: 
\|x\|^2 - \|y\|^2 \not=0\} \] 
of invertible elements has the two connected components 
\[ E^\times_\pm = \{ (x,y) \in \R^{p-1} \times \R^{q-1} \: 
\pm(\|x\|^2 - \|y\|^2) > 0\}. \] 

In $\SO_{p,q}(\R)_e$, a maximal compact subgroup is the connected 
subgroup 
\[ K = \SO_{p,q}(\R)_e \cap \OO_{p+q}(\R) \cong \SO_p(\R) \times \SO_q(\R).\] 
Hence polar decomposition shows that 
$(\SO_{p,q}(\R)_e)^h = K^h \exp(\g^h \cap \Sym_{p+q}(\R)).$ 
If an element $(k_1, k_2) \in K$ commutes with $h$, then 
it maps $e_1$ into $(\R e_1 + \R e_{p+q}) \cap (\R^p \times \{0\})$, 
hence to $\pm e_1$, and likewise with $e_{p+q}$. 
Therefore 
\[ K^h = {\rm S}(\OO_1(\R) \times \OO_{p-1}(\R)) \times 
{\rm S}(\OO_{q-1}(\R) \times \OO_1(\R)) 
\cong\OO_{p-1}(\R) \times \OO_{q-1}(\R)\] 
has four connected components. 
The group $\OO_{p-1}(\R) \times \OO_{q-1}(\R)$ acts in the obvious 
way on~$E$, preserving both connected components of $E^\times$. 
According to \cite[p.~625]{HN12}, the subgroup 
$\SO_{1,1}(\R) \cong \R^\times$ (acting on the first and last component), 
acts on $E$ by dilations. Further,
 the element $(-1, 1) \in \OO_{1,1}(\R) \subeq \OO_{p,q}(\R)$, 
 acting on the first and last component, acts on~$E$
 by inversion in the unit sphere: 
$v \mapsto \frac{v}{\beta(v,v)}$ for 
$\beta((x,y),(x,y)) = \|x\|^2 - \|y\|^2$. Since all these maps 
preserve the subsets $E^\times_\pm \subeq E^\times$, 
it follows that $G^h$ preserves both connected components~$E^\times_\pm$. 
\end{prf}

\begin{rem} \mlabel{rem:3.11} (a) The isomorphic Lie algebras
  \[ \fsl_2(\R) \cong \su_{1,1}(\C) \cong \so_{1,2}(\R)
    \cong \sp_2(\R) \cong \so^*(4) \]
  are of real rank $r = 1$, and for these we have $\g^* \cong \R$.
  These are the conformal Lie algebras of the
  $1$-dimensional Jordan algebra $\Herm_1(\K) = \R \1$
  for $\K = \R,\C,\H$. 

  \nin (b) The Lie algebra $\so_{2,2}(\R) \cong \so_{1,2}(\R)^{\oplus 2}$
  is not simple, hence not included in Table~2.
  We refer to Section~\ref{subsec:2.9} for a detailed discussion of this Lie algebra. 
\end{rem}

\subsection{Symmetric orthogonal pairs}
 \mlabel{subsec:2.3}

\begin{defn} \mlabel{def:euler-sym-pair}
 We call an orthogonal pair $(h,k)$ of Euler elements
 {\it symmetric} if $(h,k)\sim (k,h)$, i.e., there exists
$\phi \in \Inn(\g)$ exchanging $k$ and $h$ (cf.\ Definition~\ref{def:diag-act}).
\end{defn}

\begin{prop} With the  notation from {\rm Theorem~\ref{thm:1.4}},
  the orthogonal pair $(h, k^j)$ is symmetric if and only if
\begin{itemize}
\item[\rm(A)] $r$ is even and $j = r/2$.
\item[\rm(C)] $j = r$, which is the only case leading to an orthogonal pair.
\item[\rm(D)] $r$ is even and $j \in \{r-1, r\}$.
\end{itemize}
\end{prop}

\begin{prf}   
  Inspection of $\fsl_2(\R)$ shows that $(k,h) \sim (h,-k)$
  (Example~\ref{ex:sl2}), so that
 $(h,k^j)$ is symmetric if and only if, in Theorem~\ref{thm:1.4},
 the Euler element $k^j$ is conjugate under $G^h$ to $-k^j$.

 \nin {\bf Case (A):} In this case $-k^j$ is conjugate to $k^{r-j}$, so that the
 classification in Theorem~\ref{thm:1.4} shows that $(h,k^j)$ is symmetric
 if and only if $r-j = j$, i.e., if $r = 2j$.

 \nin {\bf Case (C):} In this case there is only one orbit, so that
 $(h,k^r)$ is symmetric. 

 \nin {\bf Case (D):} If $r$ is even, then the symmetry follows 
 from $-\1$ being contained in the Weyl group~$\cW(\g^*,\fc)$.
 For $r =2$, we have $D_2 \cong A_1^{\oplus 2}$, whose
   Weyl group $\cW \cong S_2 \times S_2$ also contains $-\1$.
 If $r$ is odd, then the Weyl group contains only even sign changes
 and permutations, so that $-k^r \sim k^{r-1}$, and likewise
 $-k^{r-1} \sim k^r$. So we have no symmetric orthogonal pairs in this case.
\end{prf}

\section{The abstract wedge space}
\mlabel{sec:4} 

Let $G$ be a connected Lie group
with Lie algebra $\g$ and $h \in \g$ be an Euler element.
We assume that the involution $\tau_h^\g = e^{\pi i \ad h}$ integrates
to  an involution $\tau_h$ on $G$ and form the connected Lie
group
\[ G_{\tau_h} := G \rtimes \{  \1, {\tau_h}\}.\]
We write $\eps_G \: G_{\tau_h} \to \{\pm 1\}$
for the surjective group homomorphism  with $\ker \eps_G = G$,
so that we have in particular $\eps_G(\tau_h) =  -1$.
Following \cite[Def.~2.5]{MN21}, we introduce the {\it
  abstract (Euler) wedge space of $G_{\tau_h}$} as a structure that
encodes much information about the $G$-action on standard subspaces
of Hilbert spaces, modular groups and modular conjugations.
We refer to the introduction for more details. 
Let 
\[ \cG_E := \cG_E(G_{\tau_h}) := \{ (x,\sigma)\in \cE(\g) \times G\tau_h
  \subeq \g \times G_{\tau_h}    \: \sigma^2 = e, \Ad(\sigma) = e^{\pi i \ad x}\}.\]
The group $G_{\tau_h}$ acts naturally on $\cG_E$  by 
\begin{equation}
  \label{eq:cG-act}
 g.(x,\sigma) := (\Ad^\eps(g)x, g\sigma g^{-1}),
\end{equation}
where $\Ad^\eps(g) = \eps_G(g) \Ad(g)$ is the {\it twisted adjoint action}. 
We define the {\it complementary} or {\it dual wedge} of 
$W = (x,\sigma) \in \cG_E$ by
\[ W' := (-x,\sigma) = \sigma.W.\] 
Note that $(W')' = W$ and $(gW)' = gW'$ for $g \in G$ 
by \eqref{eq:cG-act}.  
The relation $\sigma.W = W'$ is our main 
motivation to introduce  the twisted adjoint action.
This fits the relation $J_\sV \sV = \sV'$ for standard 
subspaces (\cite{Lo08}). 
It also matches the geometric interpretation in terms
of wedge regions in spacetime manifolds 
and the modular theory of operator algebras. 

We consider the pair 
\[ W_0 := (h, {\tau_h}) \in \cG_E \]
as a base point and call its  orbit 
\begin{equation}
  \label{eq:cw+b}
  \cW_+ := G.W_0 = \{ (\Ad(g)h, g {\tau_h}(g)^{-1} \cdot \tau_h) \: g \in G\}
\end{equation}
the corresponding {\it abstract wedge space}.
We write
\begin{equation}
  \label{eq:qtoOh}
  q \: \cW_+ \to \cO_h =  \Ad(G)h, \quad (x,\sigma) \mapsto x
\end{equation}
for the projection onto the first component,
\[  Z :=  Z(G), \quad Z^- := \{  z \in Z \:  \tau_h(z) = z^{-1} \},\]
and consider the homomorphism 
\begin{equation}
  \label{eq:partialW}
  \delta_h \:  G^{\{\pm h\}} := \{ g \in G \: \Ad(g)h \in \{\pm h \}\} \to Z^-,  \quad g \mapsto g \tau_h(g)^{-1}
\end{equation}
 (cf.\ \cite[\S 2.4]{MN21}). The kernel of its restriction
 to $G^h$ is the open subgroup
 \begin{equation}
   \label{eq:gw0}
   G_{W_0} = G^{h,\tau_h} = \{ g \in G \:  g.W_0 = W_0 \}
   \subeq G^h. 
 \end{equation}
 We consider the central subgroups
\begin{equation}
  \label{eq:z123}
 Z_1 := \delta_h(Z(G)) \subeq Z_2
  := \delta_h(G^h) \subeq Z_3 := \delta_h(G^{\{\pm h\}})
  \subeq Z^-.
\end{equation}
The subgroup $Z_2$ describes the orbit
$G^h.W_0 = \{h\}  \times Z_2 \cdot{\tau_h} \subeq \cW_+,$ 
which is the $q$-fiber over the Euler element~$h$. We also note that
$\cW_+ \cong G/G_{W_0} = G/G^{h,\tau_h}$.
If $h$ is not symmetric, then $G^{\{\pm h\}} = G^h$ and $Z_2 = Z_3$,
and if $h$ is symmetric, then $G^h \subeq G^{\{\pm h\}}$ is a subgroup
of index $2$. We therefore always have
\begin{equation}
  \label{eq:z3z2}
  |Z_3/Z_2| \leq 2.
\end{equation}

\begin{rem} \mlabel{rem:z2pi1oh} 
  Let $q_G \: \tilde G \to G$ denote the simply connected covering group
  (on which $\tau_h$ exists) 
and put $\tilde G_{\tau_h} = \tilde G \rtimes \{ \1,\tau_h\}$.
The subgroup $\tilde G^{\tau_h}$ is connected
(\cite[Thm.~IV.3.4]{Lo69}),  
hence equal to $\tilde G^h_e$, so that
$\tilde W_0 := (h, \tau_h) \in \cG_E(\tilde G_{\tau_h})$ has the connected
stabilizer group $\tilde G^h_e$. Therefore $\tilde \cW_+ := \tilde G.\tilde W_0$
is simply connected, and the map
\begin{equation}
  \label{eq:tildew-cov}
  \tilde q \: \tilde \cW_+ \to \cO_h, \quad (x,\sigma) \mapsto x
\end{equation}
is the simply connected covering space of the adjoint orbit $\cO_h$. We thus 
obtain the following topological interpretation of
the $Z_2$-subgroup of $\tilde G$: 
\begin{equation}
  \label{eq:pi1oh}
  \tilde Z_2 := \delta_h(\tilde G^h)
  \cong \tilde G^h/\tilde G^h_e 
  \cong \pi_0(\tilde G^h)
 \cong \pi_1(\cO_h).
\end{equation}

If $h$ is symmetric, then $-\cO_h = \cO_h$, so that we also obtain a
projective adjoint orbit
\[ \POh := \cO_h/\{ \pm \1\}, \]
for which $\tilde \cW_+$ also is the universal covering.
Accordingly,
\begin{equation}
  \label{eq:pi1z3}
  \tilde Z_3 := \delta_h(\tilde G^{\{\pm h\}}) \cong
  \pi_0(\tilde G^{\{\pm h\}}) \cong \pi_1(\POh). 
\end{equation}
\end{rem}

The complementary wedge $W_0' = (-h,{\tau_h})$
is in general not contained in $\cW_+$, but there is a replacement.
For $\alpha \in Z^-$, we call
\begin{equation}
  \label{eq:walpha'}
  W_0^{'\alpha} := (-h, \alpha {\tau_h})
\end{equation}
the {\it $\alpha$-complement of $W_0$}.

 \begin{lem}  \mlabel{lem:5.1}
  The Euler element $h$ is symmetric if and only if
  there exists an $\alpha \in Z^-$ with $W_0^{'\alpha} \in \cW_+$.
If $\Ad(g_0)h = -h$, then 
$\delta_h(g_0) Z_2 = \{ \alpha \in Z^- \: W_0^{'\alpha} \in \cW_+\}$.
\end{lem}

\begin{prf} (cf.\ \cite[Lemma~2.18]{MN21})
  For $g \in G$, the relation $g.W_0 = W_0^{'\alpha}$
  is equivalent to $\Ad(g)h = -h$ and $\delta_h(g) = \alpha$.
  This implies the first assertion. For the second
  we observe that $g_0 G^h = \{ g \in G \: \Ad(g)h = - h\}$
  and use that $\delta_h$ is a group homomorphism.
 \end{prf}

With this lemma we can also characterize when $W_0' \in \cW_+$:

\begin{lem}
  \mlabel{lem:5.2}
  The following are equivalent:
  \begin{itemize}
  \item[\rm(a)] $W_0' \in \cW_+$. 
\item[\rm(b)] There exists a $g \in G^{\tau_h}$ with $\Ad(g)h = -h$.
\item[\rm(c)] $Z_3 = Z_2$.
\end{itemize}
\end{lem}

\begin{prf} That (a) implies (b) follows from
$g.{W_0} = g.(h,{\tau_h}) = {W_0}^{'\alpha} = (-h, \delta_h(g){\tau_h}).$ 
For any $g \in G$ with $\Ad(g)h = -h$, we have
$G^{\{\pm h\}} = G^h \dot\cup g G^h$, 
so that (b) implies (c). To see that (c) implies (a), pick
$g \in G$ with $\Ad(g) h = -h$. Then (c) implies that
$\delta_h(g) = \delta_h(g_0)^{-1}$
for some $g_0 \in G^h$. Then $\Ad(gg_0)h = -h$ with
$\delta_h(gg_0) =  e$, so that ${W_0}' = gg_0.{W_0} \in \cW_+$. 
\end{prf}

\begin{rem} \mlabel{rem:z2-z3}
  (a)  (The case where $G$ is simply connected)
  If $G$ is simply connected, then $G^{\tau_h}$ is connected,
  hence coincides with $G^h_e$, so that condition (b) in
  Lemma~\ref{lem:5.2} is never satisfied. Accordingly,
  $W_0' \not\in \cW_+$ in this case, i.e., $Z_2 \not= Z_3$. 

 \nin (b)  (The case where $\g$ is simple) 
  Suppose that $\g$ is a simple Lie algebra
  and that $h$ is symmetric, so  that $\cW_+$ contains
some $\alpha$-complement $W_0^{'\alpha}$ of~$W_0$ (Lemma~\ref{lem:5.1}).
We know from \cite[Thm.~7.8]{MNO23} that $\Ad(G)^h$ is either connected
or has two connected components. In the latter case we pick
$g_- \in G^h$ with $\Ad(g_-) \not \in \Ad(G)^h_e$; otherwise
we put $g_- := e$. Then $G^h = Z(G) G^h_e \{e,g_-\}$ 
implies that
\begin{equation}
  \label{eq:z2gh}
  Z_2 = \delta_h(G^h)  = Z_1 \cup \delta_h(g_-)Z_1,
\end{equation}
so that either $Z_1 = Z_2$ (if $\Ad(G)^h$ is connected),
or $|Z_2/Z_1| = 2$.
\end{rem}

\subsection*{The central element associated to an
  orthogonal pair}

\begin{defn} \mlabel{def:2.10}

Let $G$ be a connected Lie group with Lie algebra $\g$ and
$S \subeq G$ the integral subgroup with Lie algebra
$\fs \cong \fsl_2(\R)$. 
Then $Z(S)$ is cyclic and contained in $Z(G)$
because $S$ acts on $\g$ as $\PSL_2(\R)$
  (Remark~\ref{rem:b.2}, \cite[Lemma~2.15]{MN22}), so that $Z(S)$
acts trivially on $\g$, hence is contained in $Z(G) = \ker(\Ad)$.
If $(h,k)$ is an orthogonal pair of Euler elements in
$\fs$, then $z_{h,k} := [h,k]$ is elliptic in $\fs$ and
\begin{equation}
  \label{eq:zeta-hk}
  \zeta_{h,k} := \exp(2\pi z_{h,k}) \in Z(S)
\end{equation}
is a generator of the cyclic group $Z(S)$ (see \eqref{eq:z0} in
Example~\ref{ex:sl2}). 
As $Z(S) \subeq Z(G)$, we thus attach to each 
$\Inn(\g)$-orbit of an orthogonal Euler pair $(h,k)$
the element~$\zeta_{h,k} \in Z(G)$. 
\end{defn}

\begin{rem}  Every $\fsl_2$-subalgebra $\fs \subeq \g$ 
  that contains an Euler element 
  contains a pair $(h,k)$ of orthogonal Euler elements,
  and in $\fs$ we then have the two conjugacy classes of
  the pairs $(h,k)$ and $(k,h)$ (Example~\ref{ex:sl2}).
 These pairs are distinguished by the
  $\Inn(\fs)$-invariant convex cone $C_{h,k} \subeq \fs$
  containing $[h,k]$ because   $C_{k,h} = - C_{h,k}$.
Conversely, every orthogonal Euler pair $(h,k)$ for which $(k,h)$ is
also  orthogonal generates a subalgebra $\fs_{h,k} = \Spann \{h,k,[h,k]\}$ isomorphic to $\fsl_2(\R)$ 
  (Theorem~\ref{thm:sl2-gen}) and determines a unique
  invariant cone $C_{h,k} \subeq \fs$ containing $[h,k]$.
Writing $\oEp(\g)$ for the set of orthogonal Euler pairs in $\g$
and $\SC$ for the set of all pairs
  $(\fs,C)$, consisting on a Lie subalgebra
  $\fs \cong \fsl_2(\R)$ containing an Euler element of $\g$
  and a pointed generating closed convex $\Inn(\fs)$-invariant cone
  $C \subeq \fs$,
we obtain an $\Inn(\g)$-equivariant bijection 
\[ \oEp(\g) \to \SC,
  \quad (h,k) \mapsto (\fs_{h,k}, C_{h,k}). \]
\end{rem}

\begin{ex}\mlabel{ex:2.8} (Euler elements in the Lorentz Lie algebra)
  The simple Lie algebra $\g = \so_{1,d}(\R)$, $d \geq 2$,
  is of real rank $r=1$ with $\Sigma(\g,\fa)$ of type $A_1$, so that
  there exists a {\bf single conjugacy class} of Euler elements,
  corresponding to the generators of Lorentz boosts.
  Here
  \[ \fh = \g_0(h) \cong \so_{1,1}(\R) \oplus \so_{d-1}(\R), \quad
    \g_1(h) \cong \R^{d-1}, 
    \quad \mbox{ and } \quad
     \g^*       \cong \so_{1,d-1}(\R).\]
   For $d= 2$,   $\g \cong \fsl_2(\R)$ and we have the two orbits
   of the Euler pairs  $(h_0, \pm k_0)$. Then $Z(S) = Z(G)$ is cyclic and
  $\zeta_{h_0, -k_0} = \zeta_{h_0, k_0}^{-1}$ are generators.
For $d \geq 3$, we have a non-split case with $\Sigma(\g^*,\fc)$ of
   type $C_1  = A_1$ and a single conjugacy class 
   of orthogonal pairs of Euler elements,
represented by the pair $(h,k^1)$ as
  in Theorem~\ref{thm:1.4}.  We may assume that the Lie subalgebra
  \[ \fs = \R h + \R k^1 + \R [h,k^1] \] generated by
  $h$ and $k^1$ is $\so_{1,2}(\R)$, acting on the subspace
  spanned by $\be_j$, $j = 0,1,2$.

  As $-\1 \not\in \SO_{1,d}(\R)_e$,
  $G \cong \SO_{1,d}(\R)_e$ is the centerfree group with Lie algebra
  $\so_{1,d}(\R)$, so that $K \cong \SO_d(\R)$ implies
  $Z(\tilde G) \cong \pi_1(G) \cong \pi_1(K) \cong \Z_2$.
Since $S$ contains the subgroup isomorphic to $\SO_2(\R)$,
    acting on $\R \be_1 + \R \be_2$ by rotations,
    and the natural morphism $\pi_1(\SO_2(\R)) \to \pi_1(\SO_d(\R))$
    is surjective (\cite[Prop.~17.1.10]{HN12}), 
  the homomorphism $\pi_1(S) \cong \Z \to \pi_1(G)\cong \Z_2$ is surjective.
  \end{ex}

\section{The fundamental group of the Euler orbit}
\mlabel{sec:5} 

In this section we calculate the fundamental groups
    $\pi_1(\cO_h)$ for all (not necessarily symmetric)
    Euler elements $h$ in a Lie algebra~$\g$.
    We start with the case where $\g$ is simple,
    which easily extends to the reductive case  
    (Theorem~\ref{thm:Z2-struc}).
    The general case follows by a 
    reduction procedure based on a Levi decomposition
    (Theorem~\ref{thm:leviproj}). 

\subsection{The case of simple Lie algebras} 

In this subsection we determine the fundamental group
of $\cO_h$ for all simple Lie algebras $\g$ and any 
Euler element~$h \in \g$. 
Let  $\theta$ be a Cartan involution with $h \in \fp = \g^{-\theta}$.
We consider the involution
$\tau := \tau_h \theta$ and choose $\fa \subeq \fp$ maximal
abelian with $h \in \fa$
and a Cartan subalgebra $\tilde\fc \subeq \g$ containing~$\fa$. 
Following the discussion in \cite[\S 5.1]{MNO23} and
the setting of Theorem~\ref{thm:1.4}, we find a
$(-\tau)$-invariant 
set
\[ \Gamma := \{\gamma_1, \cdots, \gamma_r\} \subeq
  \{ \beta \in \Sigma_1 := \Sigma(\g,\tilde\fc) \: \beta(h) = 1 \} \] 
of long strongly orthogonal roots. 
Let $\Gamma_0 := \Gamma^{-\tau}$
and $\Gamma_1 := \Gamma \setminus \Gamma_0$,
so that $-\tau$ acts on $\Gamma_1$ without fixed points.
For $r_0 := |\Gamma_0|$ and $r_1 := |\Gamma_1|/2$, we then have
$r = r_0 + 2r_1$ and put $s := r_0 + r_1$.
As in \cite[eq.~(24)]{MNO23}, we obtain a $\tau$-invariant subalgebra
\begin{equation}
  \label{eq:fs-decomp}
\fs  \cong \fsl_2(\R)^{r_0} \oplus \fsl_2(\C)^{r_1}
\quad \mbox{ with } \quad 
 \fs^{\tau} 
\cong \so_{1,1}(\R)^{r_0} \oplus \su_{1,1}(\C)^{r_1}.
\end{equation}
Then the subspace $\ft_\fq := \fs \cap \fq_\fk = \fs \cap \fk^{-\tau_h}$
is maximal abelian in $\fq_\fk$ and can be identified
with the compactly embedded Cartan subalgebra $\so_2(\R)^{r_0 + r_1} \subeq \fs$.
It follows from 
\cite[Prop.~5.2, Thm.~5.4]{MNO23} that 
the subalgebra $\fl$, generated by $h$ and $\ft_\fq$, is
\begin{equation}
  \label{eq:ell}
\fl =  \R h_c \oplus \fsl_2(\R)^{r_0} \oplus\fsl_2(\R)^{r_1}
\cong \R h_c \oplus \fsl_2(\R)^s \subeq \R h_c \oplus \fs,
\end{equation}
identified as in~\eqref{eq:fs-decomp}.
We also recall that $h_c = 0$ if and only if $h$ is symmetric.

Let $G := \Inn(\g) \supeq K := \Inn_\g(\fk)$, 
$L := \la \exp \fl \ra \subeq G$, and $T_Q := \exp(\ft_\fq).$
We consider the symmetric R-space 
$\cO_h^K := \Ad(K)h$ (\cite[Rem.~7.5]{MNO23}, \cite{Lo85})
and note that the inclusion
$\tilde K \into \tilde G$ implies that
\begin{equation}
  \label{eq:pi1ohk}
 \pi_1(\cO_h^K) \cong \pi_0(\tilde K^h) \cong \pi_0(\tilde G^h)
  \cong \pi_1(\cO_h) 
  \ {\buildrel \eqref{eq:pi1oh} \over \cong}\  \tilde Z_2
  = \delta_h(\tilde G^h).
\end{equation}

\begin{thm} {\rm(Structure Theorem for $\pi_1(\cO_h)$)}
  \mlabel{thm:Z2-struc}
  If $\g$ is simple, 
  then  $\pi_1(\cO_h)$ can be determined as follows: 
  \begin{itemize}
  \item $\pi_1(\cO_h)$ is trivial if $\g$ is complex or $(\g,h)$ is of
    non-split type {\rm(NST)}. 
  \item $\pi_1(\cO_h)\cong \Z$ is $\g$ is hermitian, i.e., for Cayley type {\rm(CT)}. 
  \item $\pi_1(\cO_h)\cong \Z_2$ if $(\g,h)$ is of split type {\rm(ST)}. 
  \end{itemize}
\end{thm}

\begin{prf} {\bf Step 1:} As $\pi_1(\cO_h) \cong \pi_1(\cO_h^K)$, we may use
  information on $\pi_1(\cO_h^K)$ from \cite{MNO23}.
  From \cite[Prop.~7.2]{MNO23}
it follows that the discrete subgroup
\[ \Gamma_h := \{ x \in \ft_\fq \: e^{\ad x}h = h \}
  \cong \pi_1(\cO_h^T)
  \quad \mbox{ with } \quad
  T_Q^h = \exp(\Gamma_h)\quad \mbox{ and }  \quad \cO_h^T := \Ad(T_Q)h, \]
defines a surjective homomorphism
\begin{equation}
  \label{eq:surj-homo}
\Gamma_h \cong \pi_1(\cO_h^T)
  \onto \pi_1(\cO_h^K), \quad
  x \mapsto \gamma_x, \quad
  \gamma_x(t) = e^{t \ad x}h, \ 0 \leq t \leq 1.
\end{equation}
Further $T_Q^h = Z(L')$ 
by  \cite[Lemma~7.6]{MNO23}, where
$L' \subeq L$ denotes the commutator subgroup. Writing $\fl = \R h_c \oplus
\fl_1 \oplus \cdots \oplus \fl_s$ with $\fl_j \cong \fsl_2(\R)$
as in \eqref{eq:ell}, and $L_j \subeq G$
for the corresponding integral subgroups, we thus obtain
\[ T_Q^h = Z(L_1) \cdots Z(L_s),\]
where the groups $Z(L_j)$ are cyclic. 
For $j > r_0$, the inclusion $L_j \to G$ extends by
\eqref{eq:fs-decomp} to a morphism
$\SL_2(\C) \to G$ and the corresponding adjoint
$\SU_2(\C)$-orbit of $h$ 
is diffeomorphic to $\bS^2$, hence simply connected. Therefore
$Z(L_j)$ does not contribute to $\pi_1(\cO_h^K)$
(\cite[Lemma~7.6(c)]{MNO23}).

For complex and non-split type, it thus
follows from $r_0 = 0$ that $\cO_h$ is simply connected
(cf.\ the proof of \cite[Thm.~7.8]{MNO23}).

\nin {\bf Step 2:} We are therefore left with the cases where $r_1 =0$, resp.,
$r = s$, i.e., $(\g,h)$ is of Cayley type  ($\g$ is hermitian of tube type)
or of  split type. Then the roots $\gamma_j$
are $(-\tau)$-invariant, hence identify naturally with
restricted roots in $\Sigma(\g,\fa)$.
The  set
$\Sigma_0 := \Sigma(\g,\fa) \cap h^\bot$ 
of restricted roots vanishing on $h$ is the set of
restricted roots of the Lie algebra $\g^h = \fh$.
Its Weyl group $\cW_0$ is a homomorphic image of the
normalizer of $\fa \subeq \fh_\fp$ in the connected Lie group
$\Inn_\g(\fh_\fk)$, acting on $\fh_\fp$.
In the Cayley type case the root system
$\Sigma(\g^h,\fa) \cong \Sigma(\g^*,\fc)$
is of type $A_{r-1}$, as a subsystem of 
$\Sigma(\g,\fa) \cong \Sigma(\g,\fc)$ of type $C_r$, and
in the split type case $\Sigma(\g^*,\fc)$ is of type $D_r$. 
We have he following inclusions 
\begin{align*}
 C_r
&  = \{ \pm \eps_j \pm \eps_m \: 1 \leq j,m \leq r\} \\
\supeq D_r &= \{ \pm \eps_j \pm \eps_m \: 1 \leq j\not= m \leq r\} \\
\supeq A_{r-1} &= \{ \pm (\eps_j -\eps_m) \: 1 \leq j\not= m \leq r\}
\end{align*}
with the set of long strongly orthogonal roots 
\[ \Pi = \{ 2 \eps_j \: j = 1,\ldots, r \} \subeq C_r.\]
For the root system $C_r$, the Euler element corresponds to the vector
\[ h = \frac{1}{2}(1,\ldots, 1) \in \R^r \quad \mbox{ and }\quad
  \Sigma_1 = \{ \eps_j + \eps_k \: 1 \leq j,k \leq r\}
  \subeq C_r.\] 
In both cases the Weyl group
\[ \cW^* :=  N_{\Inn(\fh_\fk)}(\fc)/Z_{\Inn(\fh_\fk)}(\fc), \]
acting on $\fc$, contains the Weyl group $S_r$ (the symmetric group), 
acting by permutations on the set $\Pi$ of strongly
orthogonal roots. We conclude that:
\begin{equation}
  \label{eq:keyfact}
  \text{ the subgroups } Z(L_j), 1 \leq j \leq r,
  \text{ are conjugate under the connected 
subgroup}\  K^h_e.   
\end{equation}  
Therefore the image of $Z(L_j)$ in $\pi_1(\cO_h^K)$
does not depend on $j$, so that the surjectivity of
\eqref{eq:surj-homo} implies that the homomorphism
$Z(L_1) \to \pi_1(\cO_h^K)$ is surjective.
The element $z_1 := [h,k_1] \in \fl_1 \cap \fk$ has the property that
$\exp(2\pi z_1)$ generates $Z(L_1)$
(cf.\ Definition~\ref{def:2.10}), so that the homotopy class of the loop
\[ \gamma \: [0,1] \to  \cO_h^K, \quad
  \gamma(t) = e^{2\pi t \ad z_1} h \]
generates $\pi_1(\cO_h^K)$. 
Accordingly,
\begin{equation}
  \label{eq:partialWz1}
  \delta_h(\exp_{\tilde G}(2\pi z_1)) = \exp_{\tilde G}(4\pi z_1)
\end{equation}
generates the subgroup $\tilde Z_2 \subeq Z(\tilde G)$.
Recall from \eqref{eq:pi1oh} that
$\tilde Z_2 \cong \pi_0(\tilde G^h) \cong \pi_1(\cO_h)$ because
$\tilde G$ is simply connected.

\nin {\bf Step 3:} {\bf Cayley type (CT):} In this case we have to show that
$\exp_{\tilde G}(4\pi z_1)$ is of infinite order. The sum
\[ z_\fk := z_1 + \cdots + z_r = [h, k_1 + \cdots + k_r] \in
  \fz(\fk) \]
(notation as in Theorem~\ref{thm:1.4}) satisfies
\begin{equation}
  \label{eq:zkz1}
  \exp_{\tilde G}(4\pi z_\fk) = \exp_{\tilde G}(r 4\pi z_1)
\end{equation}
because all central elements
$\exp_{\tilde G}(4\pi z_j), j = 1,\ldots, r$, are equal.
In fact, they are conjugate under the connected 
group $\tilde K^h_e$ by the Weyl group argument from above. 
As $Z(\tilde K)_e \cong \R$, the element \eqref{eq:zkz1} is of infinite order.

\nin {\bf Step 4:} {\bf Split type (ST):}  There are $5$ series of simple Lie algebras
of split type, containing (not necessarily symmetric) Euler elements:
\begin{equation}
  \label{eq:split-type-euler}
 \fsl_n(\R), \quad \so_{n,n}(\R), \quad \fe_6(\R), \quad \fe_7(\R), \quad 
 \mbox{ and }\quad \so_{p,q}(\R), p,q > 2
\end{equation}
(cf.\ \cite[Table~3]{MNO23}). Note that $\so_{2,q}(\R)$ is hermitian,
hence corresponds to Cayley type, and that $\so_{2,2}(\R)$ is not simple.  Therefore
we require $p,q > 2$ in \eqref{eq:split-type-euler}. 

\nin $\g = \fsl_n(\R), n > 2$: Here the conjugacy classes
of Euler elements are represented by  
\[ h = \frac{1}{p+q} \diag(q\1_p, -p\1_q), \quad p \leq q, \ n = p + q,\] 
and the strongly orthogonal roots in $\Sigma_1$ are 
$\gamma_j := \eps_j - \eps_{p+j}$, $j = 1,\ldots, p$. 
Therefore
\begin{equation}
  \label{eq:z_1insl2}
 z_1 = \frac{1}{2} \pmat{ 0 & 1 \\ -1 & 0} \in \fsl_2(\R),
  \quad \mbox{ acting\ on} \ \quad
  \R \be_1 + \R \be_{p+1}.
\end{equation}
As $\SL_n(\C)$ is simply connected,
$\ker(\eta_G) \cong \pi_1(\SL_n(\R)) \cong \Z_2$
and $\exp_{\tilde G}(4 \pi z_1)$ represents a $2\pi$-rotation in the plane
spanned by $\be_1$ and $\be_{p+1}$, hence a generator of $\pi_1(\SL_n(\R))
\subeq Z(\tilde G)$. Therefore
$\exp_{\tilde G}(4\pi z_1) \not=e$ and $e^{2 \pi \ad z_1} h = h$ imply
with \eqref{eq:partialWz1} 
that $\tilde Z_2 = \pi_1(\SL_n(\R))$, hence $\tilde Z_2 \cong \Z_2$.\\

\nin $\g = \so_{n,n}(\R), n > 2$ (for $n = 2$ it is not simple),
with restricted root system
$D_n$ and the non-symmetric Euler elements $h = h_{n-1}, h_n$,
which means that, identifying $\fa$ with $\R^n$ and
$D_n \subeq (\R^n)^*$ with the set of linear functionals of the form
$\pm \eps_j \pm \eps_k$, $j \not=k$, we have the positive system 
\[ \alpha_1 = \eps_1- \eps_2,\ \  \ldots,\ \
  \alpha_{n-1} = \eps_{n-1} - \eps_n,\ \ 
\ \alpha_n = \eps_{n-1} + \eps_n.\]
This implies that   
\[ h_{n-1} = \frac{1}{2}(1,\cdots, 1,-1) \quad \mbox{ or } \quad
  h_n = \frac{1}{2}(1,\cdots, 1).\] 
Since both Euler elements are conjugate under the full automorphism
group of $\g$, it suffices to discuss the case $h = h_n$. Then 
$\Sigma_1 = \{ \eps_j + \eps_k \: j \not=k \}$ and
\[ \gamma_1 = \eps_1+ \eps_2, \quad 
\gamma_2 = \eps_3+ \eps_4, \quad \cdots, \quad 
\gamma_{\frac{n-1}{2}} = \eps_{n-2}+ \eps_{n-1} \]
is a maximal set of strongly orthogonal roots. 
So $\fl_1 \cong \fsl_2(\R) \subeq \so_{2,2}(\R)$ acts on the linear
subspace $\Spann_\R \{ \be_1, \be_2, \be_{n+1}, \be_{n+2}\}$
(cf.\ Subsection~\ref{subsec:2.9}).

We have $\pi_1(\SO_{n,n}(\R)) \cong \Z_2^2 \subeq Z(\tilde G)$ and 
the natural homomorphism
\begin{equation}
  \label{eq:sonn-spin}
 \pi_1(\SO_{n,n}(\R)) \to \pi_1(\SO_{2n}(\C)) \cong \Z_2, \quad
 (\oline n, \oline m) \mapsto \oline n + \oline m
\end{equation}
(cf.~Remark~\ref{rem:2.19}). This shows that the kernel of
$\eta_{\tilde G} \: \tilde G \to \Spin_{2n}(\C)$ is 
\[ \ker(\eta_{\tilde G}) = \{ (\oline 0, \oline 0), (\oline 1, \oline 1) \}
  = \Delta_{\Z_2}  \subeq \pi_1(\SO_{n,n}(\R)).\]
From Remark~\ref{rem:2.19} we infer that
  $\exp(4\pi z_1) = (\oline 1, \oline 1)$, which shows that 
  $\tilde Z_2 \cong \Z_2$.\\ 

\nin $\g = \fe_6(\R)$ and the non-symmetric Euler elements
$h_1$ and $h_6$ (Theorem~\ref{thm:classif-symeuler}):
In this case Example~\ref{ex:e6-wigge} implies that $\tilde Z_2 \cong \Z_2$.\\

\nin $\g = \fe_7(\R)$:  By \cite[p.~48]{Ti67} we have
$Z(\tilde G) \cong \Z_4$, $Z(\tilde G_\C)
\cong \Z_2$ (so that $\tilde Z_1 = \tilde Z_2$),
and $\Inn(\g)^h$ is connected by
\cite[Thm.~7.8]{MNO23} because $r = 3$ is odd.
We therefore have $Z_1 = Z_2$.
From \cite[Thm.~3.1]{NO25} we know that
\[ \exp_{\tilde G_\C}(2\pi i h) \not=e,\]
so that this is a generator of $Z(\tilde G_\C) \cong \Z_2$.

From \eqref{eq:keyfact} in {\bf Step 2} above, we know that
  the exponentials  $\exp_{\tilde G_\C}(z_j)$, $j = 1,2,3$, all coincide 
so that 
\[ \exp_{\tilde G_\C}(2\pi i h) 
  = \exp_{\tilde G_\C}(2\pi(z_1 + z_2 + z_3)) = \exp_{\tilde G_\C}(6\pi z_1).\]
It follows that $\eta_{\tilde G}(\exp(6\pi z_1))$ generates $Z(\tilde G_\C)$.
Hence $\exp(6\pi z_1)$ is not of order $2$, hence
generates $Z(\tilde G) \cong \Z_4$,
and this implies that $\exp(2\pi z_1)$ generates $Z(\tilde G)$.
This in turn shows that $\exp(4\pi z_1) =  \delta_h(\exp(2\pi z_1))$
generates $\tilde Z_1 = \tilde Z_2 \cong \Z_2$.\\

\nin $\g = \so_{p+1,q+1}(\R), p,q > 1$:
Here the Euler element $h = h_{1,p+2}$ corresponds to a boost in a Minkowski
plane such as $\R\be_1 + \R \be_{p+2}$
(such as the Euler element $h = h_{1,4}\in \so_{2,2}(\R)$ in
Subsection~\ref{subsec:2.9}).
Therefore a set of strongly orthogonal roots
in $B_{p+1}$, resp., $D_{p+1}$, is given by
\[ \gamma_1 := \eps_1 + \eps_2 \quad \mbox{ and }
  \quad \gamma_2 := \eps_1 - \eps_2,\]
where $\eps_j(h) = \delta_{1,j}$. 
Then $\fl \cong \so_{2,2}(\R)$, acting on the subspace
generated by $\be_1, \be_2, \be_{p+2}, \be_{p+3}$.

As in \eqref{eq:sonn-spin} and Remark~\ref{rem:2.19}, we have
$\pi_1(\SO_{p+1,q+1}(\R)) \cong \Z_2^2$ with the natural
homomorphism to $\pi_1(\tilde G_\C) \cong \Z_2$, given by addition.
Here $z_1 \in \fl_1 \cong \fsl_2(\R)$ yields the element
$\exp(4\pi z_1) = (\oline 1, \oline 1)$
(Remark~\ref{rem:2.19}), generating $\ker(\eta_{\tilde G}) \cong \Z_2$.
We finally arrive at $\tilde Z_2 = \ker(\eta_{\tilde G}) \cong \Z_2$.
\end{prf}

  \begin{rem} Let $\g$ be a simple non-compact Lie algebra, $G := \Inn(\g)$, 
    $h \in \g$ an Euler element and $\theta$ a Cartan involution
    for which $\theta(h) = -h$. We consider the involution
    $\tau := \tau_h^\g \theta$, for which $(\g,\tau)$ is a
    non-compactly causal symmetric Lie algebra, where
    the cone $C \subeq \g^{-\tau}$ is generated by the $\Inn(\g^\tau)$-orbit of~$h$.
    We refer to \cite{MNO23} for details. 

    If $h$ is not symmetric, then $M := \cM(\g,\tau) := G/G^\tau$
    is a corresponding ``minimal'' causal symmetric space, but
    if $h$ is symmetric, then the minimal symmetric space with a
    causal structure is a $2$-fold covering
    $M := G/H$, where $H \subeq G^\tau$ is an index $2$ subgroup
    (\cite[Thm.~4.19]{MNO23}).

    For the maximal compact subgroup $K := \Inn(\g)^\theta$,
    the subgroup $K^\tau$ leaves $\fq_\fp = \g^{-\tau,-\theta}$ invariant
    and $K^\tau.h = \{h\}$ if $h$ is not symmetric, whereas
    $K^\tau.h = \{ \pm h\}$ if $h$ is symmetric
    (\cite[Lemma~4.11]{MNO23}). This implies that
    $H_K := H \cap K = K^h$ holds in both cases. 
    By \cite[Thm.~IV.3.5]{Lo69}, the map 
    $K \times \fq_\fp \to M, (k,x) \mapsto k \exp(x) H$ 
  factors through a diffeomorphism
  $K \times_{H_K} \fq_\fp \cong M$, thus defining on $M$ the structure of a
  vector bundle over $K/H_K \cong K/K^h \cong \cO_h^K$.
We conclude that, in both cases,
  \[ \pi_1(M) \cong \pi_1(K/H_K) \cong \pi_1(\cO_h^K) \cong \pi_1(\cO_h) \]
  (cf.~\eqref{eq:pi1ohk}). This means that coverings
  of $\cO_h$ are in one-to-one correspondence with coverings of
  $M$, and this has interesting consequence for locality investigations
  for local nets on $M$ and its coverings (cf.\ \cite{NO25}). 
\end{rem}

\subsection{The case of general Lie algebras}

Our results on the fundamental group $\pi_1(\cO_h)$ of the
adjoint orbit $\cO_h$ of an Euler element $h$ in a simple Lie algebra $\g$ 
immediately applies to the reductive case. Then
\[ \g = \fz(\g) \oplus \g_1 \oplus \cdots \oplus \g_n \]
with simple ideals $\g_j$. Writing 
$h = h_0 + h_1 + \cdots + h_n$ with
$h_0 \in \fz(\g)$ and $h_j \in \g_j$, the orbit
\[   \cO_h = \{h_0\} \times \cO_{h_1} \times \cdots \times \cO_{h_n} \]
is a product space. One can even go further:

\begin{thm} \mlabel{thm:leviproj} Let $h \in \g$ be an Euler element,
  $\fr \trile \g$ the solvable radical and $\fs := \g/\fr$ the
  semisimple quotient. 
  Then the restriction $q \: \cO_h \to \cO_{q(h)}$
  of the quotient map $\g \to \fs$ to $\cO_h$ is a homotopy equivalence.
  In particular $\pi_1(\cO_h) \cong \pi_1(\cO_{q(h)})$. 
\end{thm}

\begin{prf} As $\ad h$ is diagonalizable, $\g$ is a semisimple
  module for the automorphism group $\exp(\R \ad h)$.
  Hence \cite[Prop.~I.2]{KN96} implies the existence of an
  $\ad h$-invariant Levi complement $\fs \subeq \g$.
  Then $\g \cong \fr \rtimes \fs$ and $q(x,y) = y$ is the projection
  onto the second component.

  Let $G$ be a simply connected Lie group with Lie algebra $\g$ and
  $R, S \subeq G$ be the integral subgroups corresponding to $\fr$ and $\fs$,
  respectively, so that $G \cong R \rtimes S$.
  We write $h = h_r + h_s$ with $h_r \in \fr$ and $h_s \in \fs$ and observe
  that $q$ restricts to a $G$-equivariant map
  $q \:  \cO_h\to \cO_{h_s}^S,$
  where the $R$-action on $\fs \cong \g/\fr$ is trivial.
  The fiber of this map is the orbit $\Ad(G^{h_s})h$.
  We have to show that this fiber is contractible.

  As $h_s \in \fs$ and $[h,\fs] \subeq \fs$ by our choice of $\fs$, we have
  $[h_r,\fs] \subeq \fs \cap \fr = \{0\}$. So $S$ fixes $h_r$,
  and thus $S^h = S^{h_s}$. We conclude that
  $G^{h_s} = R S^{h_s} = R S^h$. 
   Therefore $\Ad(G^{h_s})h = \Ad(R)h$, which is the orbit of the
   Euler element $h$ in the solvable Lie algebra $\fb := \fr + \R h$.
   It therefore remains to show that orbits of Euler elements in
   solvable Lie algebras are contractible.

   We may therefore assume that $\g = \fr$ is solvable. We
   write $\g = \g_1 + \g_0 + \g_{-1}$ for the decomposition
   into $\ad h$-eigenspaces and $G_j = \la \exp \g_j \ra \subeq G$
   for the corresponding integral subgroups.
   Let $\fn := [\g,\g]$ and observe that $\g_{\pm 1} \subeq \fn$, so that
   $\g = \fn + \g_0$. For the corresponding normal subgroup $N \trile G$,
   we therefore have $G = N G_0$, so that $\cO_h = \Ad(G)h = \Ad(N)h$.
   The nilpotent Lie algebra $\fn$ is the sum of the two subalgebras
   $\fn_1 = \g_1$ and
   $\fn_0 + \g_{-1} = \fn_0 +  \fn_{-1}  \cong \fn_{-1} \rtimes \fn_0$.
   From \cite[Lemma~11.2.13]{HN12} it follows that
   $N = N_1 (N_0 N_{-1}) = N_1 N_0 N_{-1}$. Further
   \cite[Thm.~11.1.21]{HN12} implies that the normal subgroup $N \trile G$
   is also simply connected, so that the multiplication map induces a
   diffeomorphism
$N_1 \times (N_{-1} \rtimes N_0)\to N.$ 
   Clearly $N_0 \subeq N^h$. We claim that we have equality.
   So let $n = \exp(x_1) \exp(x_{-1})$ with $x_j \in \fn_j$ fix~$h$.
   Then
   \[ h + x_1 = h - [x_1,h] = e^{-\ad x_1} h = e^{\ad x_{-1}} h
     = h + [x_{-1},h] = h + x_{-1}\]
   leads to $x_1 = x_{-1} \in \fn_1 \cap \fn_{-1} = \{0\}$. This shows that
   $N^h = N_0$, and thus
   \[ \cO_h = \Ad(N)h \cong N/N^h \cong N_1 \times N_{-1}
     \cong \fn_1 \times \fn_{-1} \]
   is contractible.
\end{prf}

\begin{ex} (a) Let $d \geq 3$.
  If $G= \R^{1,d} \rtimes \SO_{1,d}(\R)_e = R \rtimes S$ is the
  Poincar\'e group and $h \in \so_{1,d}(\R)$ is a Lorentz boost
  (an Euler element in $\g$), and $\cO_h^S := \Ad(S)h$,
  then the $\fs$-projection yields a homotopy equivalence
  $\cO_h \to \cO_h^S$ with fiber diffeomorphic to
  $\Ad(R)h\cong \R^{d-1}$. As $(\so_{1,d}(\R),h)$ is of non-split type,
  $\tilde Z_2 \cong \pi_1(\cO_h^S)$ is trivial, so that $\cO_h$ is simply connected
  (Theorem~\ref{thm:leviproj}).

  \nin (b) This contrasts the situation in the conformal group
  $G = \SO_{2,d}(\R)$, where $\pi_1(\cO_h) \cong \Z$ follows from the
  fact that $\g = \so_{2,d}(\R)$ is simple hermitian for
  $d= 1$ or $d \geq 3$. 
\end{ex}

\section{Twisted wedge complements from orthogonal pairs}
\mlabel{sec:6}

In this section we take a closer look at twisted complements
of abstract wedges and how they are related to orthogonal
Euler pairs. Our main result is that all twists can be obtained by
  successively applying central elements $\zeta_{h,k}$
  contained in $\SL_2(\R)$-subgroups (Theorem~\ref{thm:slconj}).

\subsection{The Lie algebra $\g = \so_{2,d}(\R)$}

In this subsection we discuss the central subgroups
$Z_1 \subeq Z_2 \subeq Z_3$ for a connected Lie group~$G$ with Lie algebra
$\g =\so_{2,d}(\R)$, on which an involution $\tau_h^\g$ exists
for the Euler element $h := h_{2,3}$ that
generates a boost in the Lorentzian $(\be_2,\be_3)$-plane.

We consider the pair $W_0 := (h,\tau_h) \in \cG_E(G_{\tau_h})$.
As $h \in \so_{2,1}(\R) \cong \fsl_2(\R)$, acting on the first $3$ coordinates,
$h$ is symmetric. By Lemma~\ref{lem:5.1} there exists an 
$\alpha \in Z^-$ with $W_0^{'\alpha} \in \cW_+$.

\begin{rem} \mlabel{rem:5.1} From \cite[Thm.~7.8]{MNO23} we know that
  the group $\Ad(G)^h$ has $2$ connected components for
  $d \geq 3$ odd and that it is connected if $d$ is even.
  That it is also connected for  $d = 1,2$ 
  follows from $\so_{2,2}(\R) \cong \fsl_2(\R)^{\oplus 2}$
  and the connectedness of~$\Inn(\fsl_2(\R))^h$. 

  If $d\geq 3$ is odd, then $\Ad(G) \cong \SO_{2,d}(\R)_e$. 
  The proof of \cite[Thm.~7.8]{MNO23} shows that
  \[ g_- :=\diag(-1,-1,-1,-1,1,\ldots)
    = \exp(2\pi z_1) \]
  (cf.~\eqref{eq:ddaggz} in Remark~\ref{rem:2.19})
  is a representative of the second component of   $\Ad(G)^h$.
  Note that $g_-$ is not central in $\Ad(G)$. 
  For $\tilde g_- := \exp_{\tilde G}(2\pi z_1) \in G$, we then have
$Z_2 = Z_1 \cup \delta_h(\tilde g_-) Z_1.$ 
As $g_-^2 = e$ holds in $G$, we have $\tilde g_-^2 \in Z(\tilde G)$.
\end{rem}

From \eqref{eq:ker-etag} in the proof of Theorem~\ref{thm:4.6}, we know
that $Z_2 =\exp_G(4\pi \Z z_1)$, and from Lemma~\ref{lem:conf-gamma}  that 
\[ \exp_{\tilde G}(4\pi z_1) = (1, \oline 1)
  \in \Z \times \Z_2 \cong \pi_1(\SO_{2,d}(\R)),\]
so that
\begin{equation}
   \label{eq:dag11}
\tilde Z_2 = \{ (n,\oline n) \: n \in \Z \}
 \subeq\Z \times \Z_2= \pi_1(\SO_{2,d}(\R))
  \end{equation}
  is a subgroup of index~$2$. 
Let $r_{12}^G(\theta) \in G$ denote the unique lift of the
  $\theta$-rotations in the $(\be_1, \be_2)$-plane.
  To determine $\tilde Z_3$, we note that
  \[ \Ad(r^G_{12}(\pi)) h = r_{12}(\pi) h r_{12}(\pi) = - h,\]
  and
  \[ \delta_h (r^{\tilde G}_{12}(\pi))
    =  r^{\tilde G}_{12}(2\pi) = (1,\oline 0) \in \Z \times \Z_2.\]
  This shows that
 $\tilde Z_3 = \tilde Z_2 \cup (1,\oline 0) \tilde Z_2 = \pi_1(\SO_{2,d}(\R))$. 

For a general connected Lie group $G$ with Lie algebra $\g$,
it follows that $Z_3$ is the image of
$\pi_1(\SO_{2,d}(\R))$ in $Z(G)$.

If $d$ is odd, then $\Ad(G) \cong \SO_{2,d}(\R)_e$,
so that $Z(\tilde G) = \pi_1(\SO_{2,d}(\R)) = \tilde Z_3$,
as we have seen above. 
If $d$ is even, then $\tilde Z_3 = \pi_1(\SO_{2,d}(\R))$
is an index $2$-subgroup of 
$Z(\tilde G)$ because $Z(\SO_{2,d}(\R)_e) = \{ \pm \1\}$.
As $\Ad(G)^h$ is connected for $d$ even, we have 
$\tilde Z_1 = \tilde Z_2$ in this case.

\begin{rem}    According to Lemma~\ref{lem:5.2},
$W_0' \in \cW_+$ is equivalent to
  $Z_2 = Z_3$, i.e., to
  \begin{equation}
    \label{eq:r12tilde}
 \delta_h(r^G_{12}(\pi)) = r^G_{12}(2\pi) \in Z_2, 
\end{equation}
which, for $G \cong \tilde G/\Gamma$, is equivalent to
$r^{\tilde G}_{12}(2\pi) \in \Gamma \tilde Z_2.$
  \end{rem}

\subsection{Generating twists with $\SL_2(\R)$-subgroups}

Recall the elements $\zeta_{h,k}=\exp(2\pi z_{h,k})$, $z_{h,k} = [h,k]$
from \eqref{eq:zetahkintro}.

\begin{lem} \mlabel{lem:indep-coset}
For orthogonal pairs $(h,k)$ and $(h,k')$ of Euler elements, the
central element $\zeta_{h,k} \in Z_3$ only depends on the $G^h$-orbit of $k$, 
and
\begin{equation}
  \label{eq:zeta123}
  \zeta_{h,k,k'} := \zeta_{h,k} \zeta_{h,k'}^{-1} \in Z_2,
\end{equation}
so that the coset $\zeta_{h,k} Z_2 \subeq Z_3$ is independent of~$k$. 
\end{lem}

\begin{prf}  Recall that $\zeta_{h,k} \in Z(G)$. For $g \in G^h$, this implies
  that
  \[ \zeta_{h,\Ad(g)k} = \exp(\Ad(g) 2\pi z_{h,k}) = g \zeta_{h,k} g^{-1} = \zeta_{h,k}.\]
  This proves the first assertion. For the second, we note that  
  $\Ad(\exp(\pi z_{h,k}))h = -h$ by \eqref{eq:z0-rot},
  so that $\zeta_{h,k} = \delta_h(\exp(\pi z_{h,k})) \in Z_3.$ 
It also follows that $\exp(\pi z_{h,k}) \exp(\pi z_{h,k'})^{-1} \in G^h$, and thus 
  \[ \zeta_{h,k,k'} =
    \delta_h\big(\exp(\pi z_{h,k})\exp(-\pi z_{h,k'})\big)
    \in \delta_h(G^h) = Z_2.\qedhere\]
\end{prf}

\begin{thm} \mlabel{thm:slconj} {\rm(The $Z_3$-Theorem)}
  Let $G$ be a connected Lie group with simple Lie algebra~$\g$.
  We further assume that $(h, k^j)$, $j = 1,\ldots, m$, is a
  set of representatives of the conjugacy classes of orthogonal
  pairs of Euler elements in $\g$ as in 
  {\rm Theorem~\ref{thm:1.4}}, and that the Lie algebra involution
  $\tau_h^\g$ integrates to an involution $\tau_h$ of~$G$.
  Then the central elements 
  \[ \zeta_{h,k^j} = \exp(2\pi [h,k^j]), \quad j = 1,\ldots, m, \]
  generate the subgroup
  $Z_3 = \delta_h(G^{\{\pm h\}}) = \{ g \tau_h(g)^{-1} \: g \in G^{\{\pm h\}} \}$
and $\zeta_{h,h_1, h_2}$ generates $Z_2$. 
\end{thm}

\begin{prf} As the universal covering map
  $q_G \: \tilde G \to G$ maps the $Z_j$-groups of $\tilde G$ to those
  of $G$, it suffices to assume that $G$ is simply connected.

  Since $\zeta_{h,k^1} = \delta_h(\exp(\pi z_{h,k^1}))$ and
  $e^{\pi \ad z_{h,k^1}} h = - h$, the group
  $Z_3$ is the union of $Z_2$ and the $Z_2$-coset $\zeta_{h,k^1}Z_2$. 
  In view  Lemma~\ref{lem:indep-coset}, the coset $\zeta_{h,k^j}Z_2$
  does not depend on $j$.   Therefore it suffices to show that the elements 
 $\zeta_{h,k^j, k^m}= \exp(2\pi [h,k^j-k^m])$ generate 
  the subgroup~$Z_2$.

  We distinguish $3$ cases {\bf (A), (C), (D)},
  according to Theorem~\ref{thm:1.4}: \\
  \nin {\bf(A)} We have representatives $(h,k^j)$, $j = 0,\ldots, r$,
  and 
  \[ \zeta_{h,k^0}\zeta_{h,k^1}^{-1} = \exp(2\pi [h,k^0- k^1])
    = \exp(4\pi [h, k_1]) =  \exp(4\pi z_1)\]
  generates $Z_2$ by the proof of the Structure Theorem~\ref{thm:Z2-struc}.

  \nin {\bf(C)} All pairs are conjugate to $(h,k^r)$
  and the root system $\Sigma(\g^*,\fc)$ is of type $C_r$. 
  According to Table~2, this means that $(\g,h)$ is either complex
  or of non-split type (NST). In this case $Z_2$ is trivial
  by Theorem~\ref{thm:Z2-struc}.

  \nin {\bf(D)} For split type, with $\Sigma(\g^*,\fc)$ of type
  $D_r$, we have two representatives
  $(h,k^r)$ and $(h,k^{r-1})$ with
  \[ \zeta_{h,k^r, k^{r-1}} = \exp(2\pi [h,k^r- k^{r-1}])
    = \exp(4\pi [h, k_r]) =  \exp(4\pi z_r)\]
  and the Structure Theorem~\ref{thm:Z2-struc} implies that
  this element generates~$Z_2 \cong \Z_2$.
\end{prf}

We now extend Theorem~\ref{thm:slconj} to general Lie algebras.

\begin{cor} \mlabel{cor:z3gen} Let
  $G$ be a connected Lie group with Lie algebra $\g$,
  $h \in \g$ a symmetric Euler element for
  which the Lie algebra involution  $\tau_h^\g$ integrates to~$G$.
  Then $Z_3$ is generated by finitely many elements
  of the form $\zeta_{h,k}$, where $(h,k)$ is a pair of orthogonal
  Euler elements and $Z_2$ is generated by elements of the form $\zeta_{h,k,k'}$. 
\end{cor}

\begin{prf} The universal covering map
  $q_G \: \tilde G \to G$ maps the $Z_j$-groups of $\tilde G$ to those
  of $G$, so that it suffices to assume that $G$ is simply connected.
  
  As $h$ is symmetric, there exists an Euler element 
  $k$ for which $(h,k)$ and $(k,h)$ are orthogonal and both
    generate a Lie subalgebra $\fs_{h,k} \cong \fsl_2(\R)$
    (Theorem~\ref{thm:sl2-gen}). Then $\zeta_{h,k} \in Z_3$
  and $Z_3 = Z_2 \cup \zeta_{h,k} Z_2$. It therefore
  suffices to show that the elements $\zeta_{h,k,k'}$
  generate~$Z_2$ (cf.~Lemma~\ref{lem:indep-coset}).

Next we recall from Remark~\ref{rem:z2pi1oh} 
that $Z_2 \cong \pi_1(\cO_h)$, so that Theorem~\ref{thm:leviproj}
reduces the problem to the case where $\g$ is semisimple.
We now write $\g = \g_1 \oplus\cdots \oplus \g_L$ 
with simple ideals $\g_j$ and, accordingly,
\[ h = h_1 + \cdots + h_L \quad \mbox{ and } \quad
  k = k_1 + \cdots + k_L.\]
As $k_j \not=0$ if and only if $h_j \not=0$, it suffices to consider the
case where all $h_j$ are non-zero. Then
\begin{equation}
  \label{eq:oh-prod}
 \cO_h \cong \cO_{h_1} \times \cdots \times \cO_{h_L}
  \quad \mbox{ entails } \quad
  Z_2 \cong \pi_1(\cO_j) \cong Z_2(G_1) \times \cdots \times
  Z_2(G_L).
\end{equation}
By Theorem~\ref{thm:slconj}, there exists for each $j$
an Euler element $k_j' \in \g_j$, orthogonal to $h_j$, such that
$\zeta_{h_j, k_j, k_j'}$ generates $Z_2(G_j)$ by Theorem~\ref{thm:slconj}.
Then the Euler elements 
\[ k^j := k_1 + \cdots + k_{j-1} + k_j' + k_{j+1} + \cdots + k_L, \quad
j = 1,\ldots, L, \]
are orthogonal to $h$,  and $\zeta_{h,k,k^j}$ generates $Z_2(G_j)$, so that all these elements generate~$Z_2$ by~\eqref{eq:oh-prod}.
\end{prf}

\begin{remark}\label{rmk:gen} (Chains of twisted complements) 
Let $G$ be a connected Lie group and $h\in \cE(\fg)$ symmetric
such that $\tau_h^G \in \Aut(G)$ exists.
Using Theorem~\ref{thm:sl2-gen}, we find
$k \in \g$ such that the Lie subalgebra $\fs_{h,k}$ generated by
$h$ and $k$ is isomorphic to $\fsl_2(\R)$.
Let $z_{h,k} := [h,k]$, so that $e^{\pi \ad z_{h,k}} h = -h$
(cf.~\eqref{eq:z0-rot}). 
For $W_0=(h,\tau_h)\in \cW_+ \subeq \cG_E(G_{\tau_h})$,
we thus obtain a twisted central complement 
\[ W_0^{',\zeta_{h,k}} = \exp(\pi z_{h,k}).W_0 = (-h,\zeta_{h,k}\tau_h)
  \quad \mbox{ with } \quad
  \zeta_{h,k} \in Z_3.\]
Changing $k$ to $k'$ yields  a different twisted central complement 
$W_0^{'\alpha}$ with $\alpha = \zeta_{h,k'} = \zeta_{h,k',k} \zeta_{h,k}$.

Any twisted central complement of $W_0$ is of the form
$W_0^{'\beta}$ with $\beta \in \alpha Z_2$.
Therefore Corollary~\ref{cor:z3gen} implies the existence of a sequence
$k_1, k_1', k_2, k_2',\ldots, k_N'$ of Euler elements orthogonal to $h$
such that
\[ \beta = \zeta_{h,k_N,k_N'} \cdots \zeta_{h,k_1,k_1'} \alpha,\]
resp.,
\[ W_0^{'\beta} =
  \exp(\pi [h,k_N]) \exp(-\pi [h, k_N']) 
  \cdots   \exp(\pi [h,k_1]) \exp(-\pi [h, k_1']) W_0^{'\alpha}.\]
Geometrically, this means that we can reach any twisted central 
complement $W_0^{'\beta}$ by applying elements of suitable $\SL_2(\R)$-subgroups
to $W_0^{'\alpha}$ (cf. \cite[Thm.~3.13]{MN21} and \cite[Thm.~2.15]{MN22}). 
\end{remark}

\section{Outlook}

We have seen that the presence of a symmetric Euler element in an abstract wedge space makes the $\fsl_2(\RR)$--structure particularly relevant. In fact, the locality properties of complementary wedges can be described in terms of central elements of integral subgroups of inequivalent $\fsl_2(\RR)$--Lie algebras (Theorem~\ref{thm:slconj}). Moreover, all the
conjugacy classes of such $\fsl_2(\RR)$-subalgebras
can be classified (Theorem~\ref{thm:1.4}). Motivated by these observations, it is natural to investigate the class of generalized geometric Algebraic Quantum Field Theories (AQFTs) that can be associated with these wedge spaces. The essential link between the geometric locality properties of the Euler wedge space and the locality structure of the AQFT is provided, through the Tomita–Takesaki Theorem,
by the Bisognano–Wichmann property and 
the modular reflection property (i.e., the PCT Theorem).
This can be taken as  essential input for a consistent AQFT theory, see \cite{MN21}.

In the forthcoming work \cite{MNO25}, we prove that the Bisognano–Wichmann property imposes precise commutation relations on complementary wedge von Neumann algebras, guided by $\widetilde{\SL}_2(\RR)$-covariant subtheories. This setting also permits the formulation of a generalization of the Spin–-Statistics Theorem 
in the presence of a vacuum state, thereby providing a more general natural
point of view on the results in \cite{GL95,DLR07}.  In this context the
Spin--Statistics Theorem refers to the relation between the commutation relations of wedge algebras and the covariant representations of the symmetry Lie group.
A key role in this analysis is played by the language of standard subspaces, which provides a natural framework for formulating these structural properties. Within this broader framework, hermitian Lie algebras have a particularly intriguing role, as they contain the family of conformal symmetries of Minkowski spacetime—symmetries that arise naturally in AQFT. In our context, the role of
$\widetilde\SL_2(\RR)$ subtheories becomes particularly  transparent in describing the localization properties of nets.

Furthermore, the Bisognano--Wichmann property can be established for geometric AQFTs associated with Hermitian Lie algebras, in the spirit of \cite{BGL93}. We therefore conclude that the $\fsl_2(\RR)$ structure constitutes a fundamental building block in the analysis of  AQFTs, as developed in the present paper and its sequel \cite{MNO25}.

\appendix

\section{Some facts on unital Jordan algebras} 
\mlabel{app:1} 

We refer to \cite{FK94} for the basic concepts concerning Jordan algebras. 

\begin{defn} For a unital Jordan algebra $(E,e)$, we write 
\[ L(x) \: E \to E, \quad L(x)y = xy \] 
for its multiplication maps, $P(x) = 2 L(x^2) - L(x)^2$ 
for its {\it quadratic representation}, 
and $\detE \: E \to \R$ for the {\it Jordan determinant}. 

The group 
\[  \Str(E) := \{ g \in \GL(E) \: (\forall x \in E)\ 
P(gx) = g P(x) g^\top  \} \] 
is called the {\it structure  group of $E$}. Here 
$g^\top$ is the transpose of $g$ with respect to the 
non-degenerate symmetric bilinear form defined by 
$\beta(x,y) := \tr(L(xy))$. 
\end{defn}

Any simple unital Jordan algebra $(E,e)$ can be realized 
as $\g_1(h)$ in a simple Lie algebra~$\g$, where $h \in \g$ is a symmetric 
Euler element. There exist  
  $e \in \g_1(h)$ and $f =- \theta(e)/2\in \g_{-1}(h)$ such that
  \begin{equation}
    \label{eq:hef}
    [e, f] = h, \quad \mbox{ and also } \quad
    [h,e] = e, \quad [h,f] = -f.
  \end{equation}
On $\g_1(h)$ the bilinear product
  \begin{equation}
    \label{eq:x*y}
    x * y := [[x,f],y], 
  \end{equation}
  then defines the structure of Jordan algebra with unit~$e$ 
  (cf.~\cite{FK94}). 

The {\it conformal group} $G^* := \Conf(E)$ of  $E$ is isomorphic 
to the automorphism group $\Aut(\g)$ of~$\g$ (cf.~\cite[Prop.~2.1.3]{Be96}). 
Its identity component $\Conf(E)_e$
can be identified with the group $G := \Inn(\g) = \Aut(\g)_e$ of inner automorphisms
of~$\g$. Then $P_- := G^h \exp(\g_{-1}(h))$ is a parabolic subgroup of $G$ and
$M := G/P_-$ is a compact homogeneous space into which $E$ embeds via
\[ \eta_E \: E \to M, \quad x \mapsto \exp(x) P_-.\]
As $\eta_E(E)$ is an open dense subset of $M$, on which $G$ acts
by rational maps, we may consider $G$ as a group of rational maps on~$E$.

Then 
\[ G^*_+ := \{ g \in G^* \: (\forall x \in \dom(g))\ \dd g(x) \in 
\Str(E)_+ \} \] 
is a subgroup of $G^*_+$ because, for $g,h \in G^*_+$, the set 
$h^{-1}(\dom(g)) \cap \dom(h)= \dom(gh)$ is an open dense subset of $E$, on which we have 
\[ \dd(gh)(x) = \dd g(h(x))\dd h(x) \in \Str(E)_+. \]  
We also have 
\[ \dd g^{-1}(g(x)) = \dd g(x)^{-1} \quad \mbox{ for }\quad 
x \in \dom(g).\]

\begin{prop} \mlabel{prop:pos}
Consider the  simple real Jordan algebra $E = \g_1(h)$ 
and its conformal group $\Conf(E)\cong \Aut(\g)$, 
which acts by rational maps on $E$. 
Write $\dom(g) \subeq E$ for the domain of $g \in \Conf(E)$ in $E$. 
Then every $g \in \Conf(E)_e$ satisfies 
\[ \dd g(x)  \in \Str(E)_+ := \{ \phi \in \Str(E) \: 
\phi^*\detE \in \R_+ \detE\} \quad \mbox{ for every } \quad 
x \in \dom(g).\] 
In particular, 
\begin{equation}
  \label{eq:posdet}
 \Conf(E)_e^h \subeq \Str(E)_+.
\end{equation}
\end{prop}


\begin{prf} 
Since $E$ is a real form of the complexification of a euclidean 
Jordan algebra (\cite{BH98}, \cite{KaN19}), we obtain 
from \cite[Props.~II.3.3, III.4.2]{FK94} the relation 
\[ P(x)^*\detE  = \detE(x)^2 \detE \quad \mbox{ for } \quad x \in E,\] 
showing that $P(E^\times) \subeq \Str(E)_+$.
The inversion $j \: E^\times \to E^\times, j(x) = x^{-1}$ satisfies 
\[ \dd j(x) = - P(x)^{-1} \quad \mbox{ for } \quad x \in E^\times\] 
(cf.\ \cite[Prop.~II.3.3]{FK94}). 
Therefore $-j \in G^*_+$. Moreover, all translations 
$\tau_x(y) = x + y$, $x \in E$, are contained in $G^*_+$. 
Thus $G^*_+$ contains the translation group 
$\exp(\g_1(h))$ and also $\exp(\Ad(-j)\g_{1}(h)) = \exp(\g_{-1}(h))$. 
As the simplicity of $\g$ implies  $[\g_1(h), \g_{-1}(h)] = \g_0(h),$ 
these two subgroups generate $G \cong \Conf(E)_e$,
and we thus obtain $G \subeq G^*_+$. 
\end{prf}

\section{Pairs of Euler elements and nilpotent elements}
\mlabel{subsec:2.4}  

We briefly discuss in this subsection the relation between the classification
of orthogonal pairs of Euler elements (Theorem~\ref{thm:1.4})
and the well-known 
classification of nilpotent elements (\cite[Ch.~8]{Bo90b}). 

Let $\g$ be a simple non-compact real Lie algebra. 
In \cite[Ch.~VIII, \S 11]{Bo90b}, $\fsl_2$-triples are
defined as triple $(x,h,y)$ with $[h,x] = 2x, [h,y] = 2y$ and
$[x,y] = -h$. We write $\fs(x,h,y) := \Spann_\R \{x,h,y\}$
for the corresponding subalgebra. 

If $(h,k)$ is an orthogonal pair of Euler elements, then we
write $k \in \g_1(h) + \g_{-1}(h)$ accordingly
as $k = k_1 + k_{-1}$. Then $k_1 \in \g_1(h)$ is nilpotent and
$(k_1, 2h)$ can be enlarged to an $\fsl_2$-triple
which is a basis of $\Spann_\R \{ h,k,[h,k]\} \cong \fsl_2(\R)$
(\cite[Thm.~3.13]{MN21}).

If, conversely, $x \in \g_1(h)$, then $x$ is a nilpotent element,  
and if $h \in [x,\g]$, resp., $h \in [x,\g_{-1}(h)]$, then 
there exists a unique $y \in \g_{-1}(h)$ with $[x,y] = 2h$, so that
$(x,2h,-y)$ is an $\fsl_2$-triple 
(Jacobson--Morozov Theorem, \cite[Ch.~VIII, \S 11, no.~2, Prop.~2]{Bo90b}). 
Then 
\[ k := \frac{x  + y}{2} \]
satisfies
\[ [k,h] = \frac{1}{2}(y - x), \quad \Big[k,\frac{y-x}{2}\Big]
  = \frac{1}{4} [x+y,-x+y] =  \frac{1}{2}[x,y] = h, \]
which implies that $k \in \fs(x,2h,-y)$ is an Euler element, hence 
conjugate to $h$ (Example~\ref{ex:sl2}) 
and therefore also an Euler element in $\g$. 
By construction, $\tau_h^\g(k) = - k$, so that $(h,k)$ is orthogonal.

For a fixed Euler element $h$, we
thus obtain a bijection between elements $x \in \g_1(h)$
with $h \in [x,\g_{-1}(h)]$ and orthogonal Euler pairs $(h,k)$. 
As this bijection is $G^h$-equivariant, the classification of
orthogonal pairs of Euler elements corresponds to the
classification of $G^h$-orbits in
\begin{equation}
  \label{eq:e-times1}
  E^\times  := \{ x \in \g_1(h) \: h \in [\g_{-1}(h),x]\}.
\end{equation}
If $x \in E^\times$, then there exists an $\fsl_2$-triple 
$(x,2h,-y)$, and the representation theory of $\fsl_2(\R)$ implies that
$\ker(\ad x)^2 \cap \g_{-1}(h) = \{0\}$. This condition
is equivalent to the map
\[ (\ad x)^2 \: \g_{-1}(h) \to \g_1(h) \]
being injective, hence bijective, as $\dim \g_1(h) = \dim \g_{-1}(h)$.
Therefore
\begin{equation}
  \label{eq:e-times2}
 E^\times =   \{ x \in \g_1(h) \:  (\ad x)^2 \g_{-1}(h) = \g_1(h)\}
 =     \{ x \in \g_1(h) \:  \ker((\ad x)^2) = \g_0(h) + \g_1(h)\}
\end{equation}
is the set of invertible elements in the Jordan algebra $E$
(cf.\ \cite[\S 5.3]{BN04}).

We thus obtain a bijection between $G$-orbits of orthogonal
Euler pairs and $G^h$-orbits in the open subset $E^\times$ of invertible
elements in the Jordan algebra $E$.

\begin{rem}
As the above correspondence maps $-x \in E^\times$ to $-k$, an orthogonal
Euler couple $(h,k)$ is symmetric if and only if the invertible
element $x \in E$ is symmetric in the sense that $-x \in \Ad(G^h)x$.
\end{rem}

\begin{rem} \mlabel{rem:b.2} For $h, x \in \g$ with $[h,x] = x$
  and an $\fsl_2$-triple $(x,2h,y)$, 
  the condition $\Spec(\ad h) \subeq \{ \pm 1,0\}$ is equivalent to
$h$ being an Euler element in $\g$, and then
\[  \Inn_\g(\fs(x,2h,y)) \cong \PSL_2(\R) \quad
  \mbox{ and }  \quad (\ad x)^3 = 0.\] 
Conversely, these two conditions imply that 
$\Spec(\ad h) \subeq \{\pm 1,0\}$, so that $h$ is an Euler element.
\end{rem}

\begin{rem} \mlabel{rem:b3} (a) If $\g$ is simple and $h \in \g$ a symmetric Euler
  element, then
we have seen in the proof of Theorem~\ref{thm:1.4} that
    $E := \g_1(h)$ carries a natural unital Jordan algebra
    structure, whose unit $e$ generates, together with $\theta(e)$, a
    Lie subalgebra isomorphic to $\fsl_2(\R)$. The proof
    of Theorem~\ref{thm:1.4} also 
    shows that the orbits of Euler pairs $(h,h_j)$ correspond to
    the connected components of the open subset $E^\times$ of invertible
    elements in~$E$ (cf.\ \eqref{eq:e-times2}).
    
    If $\Sigma(\g^*,\fc)$ is of type $A_{r-1}$ 
    (CT), then $E$ is a euclidean Jordan algebra or rank~$r$ and
    $E^\times$ has $r+1$ connected components, corresponding to
    elements of signature $(p,q)$, $p + q = r, p,q \in \N_0$.
    For root systems of type $D_k$ (ST), the set $E^\times$ has $2$ connected
    components, and, for type $C_r$, it is connected, which happens
    for (NST) and complex type. 

\nin (b) The classification of the $G^h$ orbits
  in the open subset $E^\times$ of invertible elements in $E$ 
 matches the double cosets of the Weyl group pair 
$\cW^* := \cW(\g^*,\fc) \subeq \cW := \cW(\g,\fc)$, arising from the 
pairs $(A_{r-1}, C_r)$, $(C_r, C_r)$ and $(D_r, C_r)$. 
For type (A) $\cW^* \cong S_r \subeq \cW \cong \Z_2^r \rtimes S_r$, 
for type (C) $\cW^* =  \cW$, and for type (D) 
$\cW^*\subeq \cW$ is an index $2$ subgroup. Therefore the orbits
of pairs of orthogonal Euler elements 
corresponds to the double cosets in $\cW^*\backslash \cW/\cW^*$.
\end{rem}

\section{The conformal Lie algebra $\so_{2,d}(\R)$}
\mlabel{subsec:2.9}

In this appendix we discuss the Lie algebras
  $\so_{2,d}(\R)$ in some detail. This serves two purposes.
  It is the conformal Lie algebra of $d$-dimensional Minkowski space,
  hence relevant for applications in Physics,
  and we  also need some concrete information on
  the generators of its fundamental group
  in Section~\ref{sec:5}. 
  
The simple Lie algebra $\g = \so_{2,d}(\R)$, $d \geq 3$, 
  is of real rank $r=2$ with $\Sigma(\g,\fa)$ of type $C_2$, so that 
  there exists a single conjugacy class of Euler elements
  (Theorem~\ref{thm:1.4}).
  The Lie algebra $\g$ is hermitian  of tube type
  and $(\g, \tau_h^\g)$ is of Cayley type.   Here
  \[ \fh = \g_0(h) \cong \so_{1,1}(\R) \oplus \so_{1,d-1}(\R), \quad
    \g_1(h) \cong \R^{1,d-1}, 
    \quad \mbox{ and } \quad
     \g^*       \cong \so_{1,d-1}(\R).\]
    As $\Sigma(\g^*,\fc)$ of type $A_1$, we obtain $3$ orbits
    of orthogonal pairs $(h,k^j)$, $j = 0,1,2$, corresponding
    to the three open orbits of the connected Lorentz group
    $\SO_{1,d-1}(\R)_e$ in Minkowski space:
   the positive/negative timelike vectors and the spacelike vectors
 (see also Subsection~\ref{subsec:2.4} for a generalization of this picture).
The Euler elements     from \eqref{eq:hj} satisfy 
    $k^2 = - k^0$ and $(h,k^1) \sim (h,-k^1)$. 

    Consider the subalgebra $\fs = \fs_1 \oplus \fs_2
    \cong \fsl_2(\R)^{\oplus 2}$ and recall $h_0$ and $k_0$ from \eqref{eq:h0k0}.
      Then $h = (h_0, h_0),$ and 
      \[       k^2 = (k_0, k_0) = k_1 + k_2 = -  k^0 \quad \mbox{ and } \quad 
    k^1 = (k_0, -k_0) = k_1 - k_2 
 \]
in the notation from Theorem~\ref{thm:1.4}.
 We describe matrix models of $k_{1/2}$ below.
The Lie subalgebra $\fs$ intersects the closed convex invariant cone $C_\g$
  in $\so_{2,d}(\R)$ in 
  \[ C = C_0 \oplus C_0 \quad  \mbox{ with } \quad C_0^\circ
    \ni z_0 = [h_0, k_0] = \frac{1}{2} \pmat{ 0 & 1 \\ -1 & 0}\]
  (cf.\ Example~\ref{ex:sl2}). 
  Then we have 
  \begin{align*}
[h, k^2] &= ([h_0, k_0], [h_0, k_0]) = (z_0,z_0) \in C^\circ, \\ 
[h, k^1] &= ([h_0, k_0], -[h_0, k_0]) = (z_0,-z_0) \not\in \pm C, \\ 
[h, k^0] &= - [h,k^2] \in -C^\circ.
  \end{align*}
  Let $\fs^j \subeq \fs$ be the Lie subalgebra generated by $h$ and $k^j$, and
  $S^j \subeq \tilde G$ the integral subgroup  generated by
  $\exp(\fs^j)$.
  We now determine the central elements 
$\zeta_{h,k^j}$, $j = 0,1,2.$ 
  By construction, we have $\zeta_{h,k^2} = \zeta_{h,k^0}^{-1}$,
  so that it suffices to consider the cases $j = 0$ and $j = 1$.

  \begin{lem} In $\so_{2,2}(\R)$ the two commuting boost Euler elements
    $h := h_{1,4}$ and $k := h_{2,3}$ (cf.\ \eqref{eq:hjk}) 
  have the property $h$ and $\shalf(h\pm k)$ represent all 
  conjugacy classes of Euler elements. 
The orthogonal pairs containing $h$ are represented by $(h,k^j)$ with 
  \[ k^2 =- k^0 = h_{1,3} \quad \mbox{ and } \quad
    k^1 = - h_{2,4}.\]
The corresponding central elements of $\tilde\SO_{2,2}(\R)_e$ are 
  \[ \zeta_{h,k^1} = (1,0) \quad \mbox{ and } \quad
    \zeta_{h,k^2}  = \zeta_{h,k^0}^{-1} = (0,1) \]
  in
  \[     \Z^2 \cong \pi_1(\SO_2(\R)) \times \pi_1(\SO_2(\R))
    \cong \pi_1(\SO_{2,2}(\R)) \subeq Z(\tilde\SO_{2,2}(\R)_e).\]
\end{lem}

\begin{prf} First we observe that
  $\g_0(h) = \R h \oplus \R k = \fa$.
  From
  \[ \ad h \pmat{
      0 & a & b_{11} & b_{12} \\
      -a & 0 & b_{21} & b_{22} \\
      b_{11} & b_{21} & 0 & d \\
      b_{12} & b_{22} & -d & 0}
    =  \pmat{
      0 & b_{22}  & -d & 0 \\ 
      -b_{22} & 0 & 0 & a \\
      -d & 0 & 0 & -b_{11} \\
      0 & a & b_{11} & 0}\] 
  we derive that
  \[ \g_1(h) = \R e_1 + \R e_2 \quad \mbox{ for } \quad
    e_1 = \frac{1}{4}\pmat{
      0 & -1 & 1 & 0\\
      1 & 0 & 0 & -1\\
      1 & 0 & 0 & -1\\
      0 & -1 & 1 & 0}, \quad 
    e_2 = \frac{1}{4}\pmat{
      0 & 1 & 1 & 0\\
      -1 & 0 & 0 & 1\\
      1 & 0 & 0 & -1\\
      0 & 1 & 1 & 0}.\] 
  We further have $[k,e_1] = e_1$ and $[k,e_2] = -e_2,$ 
  so that we obtain the root spaces $\g_{\gamma_j}$: 
  \[ \g_{\gamma_1} = \R e_1, \quad \gamma_1(h) = \gamma_1(k) = 1
    \quad \mbox{ and } \quad 
    \g_{\gamma_2} = \R e_2, \quad \gamma_2(h) = -\gamma_2(k) = 1.\]
This leads to the $\fsl_2$-subalgebras
\[ \fs_1 = \g_{\gamma_1} + \R (h+k)  + \g_{-\gamma_1} \quad \mbox{ and }  \quad 
  \fs_2 = \g_{\gamma_2} + \R (h-k)  + \g_{-\gamma_2}.\]
  It follows in particular that $\frac{1}{2}(h \pm k)$ are
  are also Euler elements, but with $4$-dimensional centralizer. Therefore 
  any Euler element in the semisimple Lie algebra
  $\so_{2,2}(\R) \cong \fs_1 \oplus \fs_2 \cong \fsl_2(\R)^{\oplus 2}$
  is conjugate either to $h$ (if its centralizer is abelian),
  or to $\frac{1}{2}(h\pm k)$ (Euler elements in the ideals
  $\fs_1$ and $\fs_2$). 
We further obtain 
  \[ k_1 = e_1 + e_1^\top
= \frac{1}{2} \pmat{
      0 & 0 & 1 & 0\\
      0 & 0 & 0 & -1\\
      1 & 0 & 0 & 0\\
      0 & -1 & 0 & 0}, \quad 
    k_2 = e_2 + e_2^\top = \frac{1}{2}\pmat{
      0 & 0 & 1 & 0\\
      0 & 0 & 0 & 1\\
      1 & 0 & 0 & 0\\
      0 & 1 & 0 & 0},\]
  and since these symmetric matrices are conjugate in
  $\so_{2,2}(\R)$ to $\frac{1}{2}(h \pm k)$, they are both
  Euler elements. 
  We thus obtain the representatives $(h,k^j)$ of orthogonal pairs of
  Euler elements: 
  \[ k^2 = k_1 + k_2 = h_{1,3}, \quad k^1 = k_1 - k_2 = -h_{2,4}
    \quad \mbox{ and } \quad
    k^0 = -k_1 - k_2 = - k^2.\] 
  This leads to
  \begin{equation}
    \label{eq:hxbrack1}
z_{h,k^2} =  [h, k^2] = [h, e_1 + e_1^\top + e_2 + e_2^\top]
    = e_1 - e_1^\top + e_2 - e_2^\top
   =  \pmat{
      0 & 0 & 0 & 0\\
      0 & 0 & 0 & 0\\
      0 & 0 & 0 & -1\\
      0 & 0 & 1 & 0} 
  \end{equation}
and
\begin{equation}
  \label{eq:hxbrack2}
z_{h,k^1} =   [h, k^1] = [h, e_1 + e_1^\top - e_2 - e_2^\top]
    = e_1 - e_1^\top - e_2 + e_2^\top
   =  \pmat{      0 & -1 & 0 & 0\\
      1 & 0 & 0 & 0\\
      0 & 0 & 0 & 0\\
      0 & 0 & 0 & 0}. 
\end{equation}
These elements generate rotations in the $\be_1,\be_2$, and the
$\be_3, \be_4$-plane, so that
$\zeta_{h,k^1}, \zeta_{h,k^2} 
\in \pi_1(\SO_{2,2}(\R)) \subeq Z(\tilde\SO_{2,2}(\R)_e)$ 
correspond to generators of the subgroups $\pi_1(\SO_2(\R)) \cong \Z$.
\end{prf}

For {\bf $d > 2$}, it remains to describe the homomorphism
  $\pi_1(\SO_{2,2}(\R)) \to \pi_1(\SO_{2,d}(\R))$.   
  In the Minkowski space $E \cong \R^{1,d-1} \cong \g_1(h)$,
  the two elements $e_j \in \g_{\gamma_j}$, $j = 1,2$, correspond
  to positive light rays, resp., a Jordan frame for the euclidean Jordan algebra
  structure. The Lie algebra generated by $e_1, e_2, \theta(e_1)$ and
  $\theta(e_2)$ coincides with $\fs \cong \conf(\R^{1,1}) \cong \so_{2,2}(\R)$. 

  \begin{lem} Assume that $d \geq 3$. With the canonical
    identifications $\pi_1(\SO_{2,2}(\R)) \cong \Z^2$ 
    and 
    \[ \pi_1(\SO_{2,d}(\R)) \cong \pi_1(\SO_2(\R)) \times \pi_1(\SO_d(\R))
      \cong \Z \times \Z_2, \]
    the natural homomorphism
    $ \pi_1(\SO_{2,2}(\R)) \to \pi_1(\SO_{2,d}(\R))$ takes the form
  \begin{equation}
    \label{eq:homo22c}
    (n,m) \mapsto (n, \oline m), \quad \Z^2 \to \Z \times \Z_2,
  \end{equation}
where $\oline m = m + 2 \Z$ denotes the congruence class of $m$ modulo~$2$.
  \end{lem}

  \begin{prf} With the stated identifications, and the observation
    that $\SO_2(\R) \times \SO_d(\R)$ is a maximal compact subgroup
    of $\SO_{2,d}(\R)_e$ (as a consequence of the polar decomposition;
    see \cite[Prop.~4.3.3]{HN12}),
    the assertion boils down to the well-known fact that
    the inclusion     $\SO_2(\R) \into \SO_d(\R)$
    induces for $d \geq 3$ the unique surjective homomorphism
    \[ \Z \cong \pi_1(\SO_2(\R)) \onto \pi_1(\SO_d(\R)) \cong \Z_2 \]
    (\cite[Prop.~17.1.10]{HN12}). 
  \end{prf}

     In these terms, we obtain from
  \eqref{eq:hxbrack1} and \eqref{eq:hxbrack2}: 
  \begin{lem} \mlabel{lem:conf-gamma} For $d \geq 3$, we have 
    \begin{equation}
    \label{eq:gammasb}
    \zeta_{h,k^0 } =  \zeta_{h,k^2 } = (0, \oline 1)
    \quad \mbox{ and } \quad
    \zeta_{h,k^1 } = (1, \oline 0)
\quad \mbox{        in }\quad  \Z \times \Z_2 \cong \pi_1(\SO_{2,d}(\R)), 
  \end{equation}
and these elements generate the group~$\pi_1(\SO_{2,d}(\R))$. 
\end{lem}

\begin{rem} \mlabel{rem:2.19}
For $p,q \geq 3$, we have the natural homomorphism 
\[ \pi_1(\SO_{p,q}(\R)) \cong \Z_2 \times \Z_2 \to \pi_1(\SO_{p+q}(\C)) \cong
  \Z_2, \quad (\oline n, \oline m) \mapsto \oline n + \oline m,\]
and, for $p = 2$, $q \geq 3$, 
\[ \pi_1(\SO_{2,q}(\R)) \cong \Z \times \Z_2 \to \pi_1(\SO_{2+q}(\C)) \cong
  \Z_2, \quad (n, \oline m) \mapsto \oline n + \oline m.\]
For
  $G = \tilde\SO_{p,q}(\R)_e$ and its universal complexification 
  $\eta_G \: G \to \Spin_{p+q}(\C)$, this shows that 
$\ker(\eta_G) = \{ (\oline n, \oline n) \: n \in \Z \} \cong \Z_2$ 
in the first case, and, for $p = 2$, we have
\[ \ker(\eta_G) = \{ (n, \oline n) \: n \in \Z \} \cong \Z. \] 
This group is always cyclic, generated by~$(1,\oline 1)$ for $p = 2$,
and by $(\oline 1, \oline 1)$ for $p > 2$. 

On the other hand, we have in the notation of Theorem~\ref{thm:1.4}
in $\so_{2,2}(\R)$ 
\begin{equation}
  \label{eq:ddaggz}
 z_1 := [h, k_1] 
  = \frac{1}{2}[h, k^2 + k^1]
  = e_1 - e_1^\top
  =  \frac{1}{2}\pmat{
      0 & -1 & 0 & 0\\
      1 & 0 & 0 & 0\\
      0 & 0 & 0 & -1\\
      0 & 0 & 1 & 0},
\end{equation}
  so that, for $p = 2$,  $\exp_G(4\pi z_1)$ represents
  the element $(1,\oline 1) \in \Z \times \Z_2 \cong
  \pi_1(\SO_{2,d}(\R))$, which generates $\ker(\eta_G)$.
\end{rem}

\section{$\pi_1(\cO_h)$ by
  Wiggerman's Theorem}
\mlabel{app:c}

The following is a special case of Wiggerman's Theorem \cite[Thm.~1.1]{Wi98}
on fundamental groups of flag manifolds. 
\begin{thm}{\rm(Wiggerman's Theorem)} \mlabel{thm:wiggeman}
  Let $\g$ be a simple real Lie algebra and
  $G = \Inn(\g)$. We assume that 
  $S = \{ \alpha_1, \ldots, \alpha_r\}\subeq \Sigma(\g,\fa)$
  is a set of simple roots and that $x_1, \ldots, x_r \in \fa$ are the dual
  basis, i.e., $\alpha_j(x_i) = \delta_{ij}$.   Put
  \[  S^* := \{ \alpha_j \: \dim \g_{\alpha_j} = 1\}.\] 
  Let $P_j \subeq G$ be the parabolic subgroup with Lie algebra
$\fp_j = \g_0 + \sum_{\alpha(x_j) \leq 0} \g_\alpha.$     Then
\[ \pi_1(G/P_j) \cong
      \begin{cases}
        \Z & \text{ if }  \alpha_j \in S^* \ \text{ and } \
             (\forall \alpha_i \in S^*)\  \alpha_j(\alpha_i^\vee) \in 2 \Z\\
        \Z_2 & \text{ if }  \alpha_j \in S^* \ \text{ and } \
             (\exists \alpha_i \in S^*)\  \alpha_j(\alpha_i^\vee) \in 2 \Z+1\\
        \{e\} & \text{ if }  \alpha_j \not\in S^*.
      \end{cases}\] 
  \end{thm}

  \begin{prf} This follows immediately from Wiggerman's result
    (Theorem~\cite[Thm.~1.1]{Wi98}) 
    that deals with more general flag manifolds. He describes the
    fundamental group of $G/P_j$ in terms of generators
    $t_\alpha$, $\alpha \in S^*$, with the relations that in our case
    boil down to $t_\alpha = e$ for $\alpha \not= \alpha_j$ and
    \[ t_\beta t_\alpha = t_\alpha t_\beta^{\eps(\alpha,\beta)}, \quad
      \eps(\alpha,\beta) := (-1)^{\beta(\alpha^\vee)}
      \quad \mbox{ for } \quad \alpha \not=\beta.\]
    Applying this to $\beta = \alpha_j$ and $ \alpha = \alpha_i$,
    it follows that the fundamental group is generated by~$t_{\alpha_j}$
    and that this element satisfies the non-trivial relation $t_{\alpha_j}^2 = e$
    if and only if~$S^*$ contains an element
    $\alpha_i$ with $\alpha_j(\alpha_i^\vee)$ odd.
\end{prf}

  \begin{rem} \mlabel{rem:wigge} (a) If $\g$ is a split real form of $\g_\C$,
    then all root spaces $\g_\alpha$ are $1$-dimensional, so that
    $S = S^*$.     

    \nin (b) If $\Sigma(\g,\fa)$ is simply laced, i.e.,
    of type $ADE$, and $r > 1$,
    then there exists a simple root $\alpha_i$ corresponding to a neighbor
    of $\alpha_j$ in the Dynkin diagram. We then have
    $\alpha_j(\alpha_i^\vee) =-1$, so that we obtain
    $\pi_1(G/P_j) \cong \Z_2$.

    \nin (c) If $h \in \fa$ is an Euler element,
    then $h = x_j$ for some $j$ and $\alpha_j$ is the unique
    simple root with $\alpha_j(h) = 1$. In this case
    \[ P_j \cong G^h \exp(\g_{-1}(h))\]
    is homotopy equivalent to $G^h$, so that
    \begin{equation}
      \label{eq:pi1ohb}
      \pi_1(\cO_h) \cong \pi_1(G/P_j).
    \end{equation}

    \nin (d) To see the generator of $\pi_1(G/P_j)$ in terms of
    Lie subalgebras, one may consider the subgroup
    $G_j \subeq G$ whose Lie algebra is generated by
    the root spaces $\g_{\pm \alpha_j}$. We only have a non-trivial
    generator if these root spaces are $1$-dimensional, 
    so that $\g_j \cong \fsl_2(\R)$.
    Then the $G_j$-orbit $C_j\cong \bS^1$
    of the base point in $G/P_j$ is a circle
    whose homotopy class generates $\pi_1(G/P_j)$
    (see \cite{Wi98} for details).

    If $P_j = U_j \rtimes L_j$ is its Levi decomposition, then
    $L_j = G^{x_j}$ with $\alpha_i(x_j) = \delta_{ij}$.
From 
\[ \fa = \ker \alpha_j \oplus \R x_j
  = \ker \alpha_j \oplus \R \alpha_j^\vee, \quad
  \alpha_j^\vee \in \g_j \quad \mbox{ and } \quad
      [\g_j, \ker \alpha_j] = \{0\},\] it
    follows that $G_j^{x_j} = G_j^{\alpha_j^\vee}$. Therefore
    the $G_j$-orbit $C_j \subeq G/P_j$ can be identified with
    the corresponding flag manifold of $G_j$. In particular we
    have a short exact sequence 
    \[\underbrace{\pi_1(G_j)}_{\cong \Z} \into
      \underbrace{\pi_1(C_j)}_{\cong \Z} \onto
      \underbrace{\pi_0(G_j^{x_j})}_{\cong \Z_2 \ \text{or}\ \{e\}}.\]
If $G_j \cong \PSL_2(\R)$, then $G_j^{x_j}$ is connected, and if       
$G_j \cong \SL_2(\R)$, then $G_j^{x_j}$ has $2$ components.
Therefore $\pi_1(G_j) \to \pi_1(C_j)$ is surjective if and only if
$G_j \cong \PSL_2(\R)$. We conclude that
\[ \pi_1(G_j) \to \pi_1(G/P_j) \]
is surjective if and only if $G_j \cong \PSL_2(\R).$

\nin (e) In the context of Euler elements $h = x_j$
and $G_j = L_1$ as above,
where $\gamma_1$ is conjugate to a simple root, it follows that
$\pi_1(L_1) \to \pi_1(\cO_h)$ is surjective if and only if $L_1 \cong \PSL_2(\R)$.
This in turn is equivalent to
\[ \Sigma(\g,\fa)(\gamma_1^\vee) \subeq 2 \Z.\]
If $r \geq 2$, then $\Sigma(\g,\fc)$ is of type
$C_r$ or $BC_r$, so that there exists a root $\beta$ with
$\beta(\gamma_1^\vee) = 1$. This implies that $L_1 \cong \SL_2(\R)$,
and therefore the homomorphism
$\pi_1(L_1)\to \pi_1(\cO_h) \cong \Z_2$ vanishes.
\end{rem}

 \begin{ex} \mlabel{ex:e6-wigge}
    Consider  the split real form $\g = \fe_6(\R)$ of type $E_6$
    and the corresponding simply connected Lie group $G$. 
We obtain with Remark~\ref{rem:wigge}(b) that
    $\pi_1(\cO_h) \cong \Z_2$ for all Euler elements $h\in \g$.     
    With \eqref{eq:pi1ohk} we thus obtain
    $Z_2   \cong \Z_2$.
    Moreover,  $Z(G) \cong \Z_2$ (\cite[p.~46]{Ti67}) implies
    that $\tau_h$ acts trivially on $Z(G)$, so that
    $ Z_1 = \{e\}$. In view of Remark~\ref{rem:z2-z3}(b),
    $Z_2 \not = Z_1 = \{e\}$, and therefore  $\Inn(\g)^h $
    is {\bf not connected} for $\g = \fe_6(\R)$.\begin{footnote}
{This case is missing in \cite[Thm.~7.8]{MNO23}.}
    \end{footnote}
  \end{ex}

\begin{rem} It is instructive to see how the fundamental 
  group $\pi_1(\cO_h)$ can be determined by
  Wiggerman's Theorem~\ref{thm:wiggeman}.
  It suggests the following algorithm. We start with an Euler
  elements $h_j$ as in Theorem~\ref{thm:classif-symeuler}.
  Then there exists a unique simple root
  $\alpha_j$ with $\alpha_j(h_j) = 1$, and all other simple roots vanish
  on $h_j$. We have the following possibilities:
  \begin{itemize}
  \item   If $m_{\alpha_j} = \dim \g_{\alpha_j} > 1$, then the
  fundamental group is trivial. 
\item   If $m_{\alpha_j} = 1$ and there exists another simple
  root $\alpha_i$ with $m_{\alpha_i} = 1$ and
  $\alpha_j(\alpha_i^\vee)$ odd, then   $\pi_1(\cO_h) \cong \Z_2$.
\item   If $m_{\alpha_j} = 1$ and all other simple 
  roots $\alpha_i$ either satisfy $m_{\alpha_i} > 1$ or
  $\alpha_j(\alpha_i^\vee)$ even, then $\pi_1(\cO_h) \cong \Z$.
  \end{itemize}

  \nin {\bf Complex type:} If $\g$ is a complex simple Lie algebra,
  then $m_{\alpha} = 2$ for all roots $\alpha$, so that $\cO_h$ is
  simply connected. \\

  \nin {\bf Cayley type:} If $\g$ is hermitian, then the multiplicity
  of the strongly orthogonal restricted roots is $1$
  (cf.\ \cite[Prop.~3.3]{MNO23}),
  and for the other restricted roots (and $\g \not=\so_{2,d}(\R)$)
  it is contained in $\{1,2,4,8\}$,
  which is the dimension of the coordinate field
  $\R,\C,\H$ or $\bO$ of the corresponding euclidean simple
  Jordan algebra. For $\g = \so_{2,d}(\R)$, the multiplicity is~$d-2$.

  For $\g = \so_{2,1}(\R) \cong \fsl_2(\R)$,  the root system if of
  rank~$1$, so that $r = 1$ and $\pi_1(\cO_h)$ is infinite cyclic.
  The only other simple Lie algebra for which two simple roots of multiplicity~$1$
  exist is $\g =  \sp_{2n}(\R)$, corresponding to the Jordan algebra
  $\Sym_n(\R)$. In this case we have to show that $\alpha_j(\alpha_i^\vee)$
  is always even. Here the system of restricted roots is of type
  \[ \Sigma(\g,\fa) \cong C_r
    = \{ \pm \eps_i \pm \eps_k \: 1 \leq i,k \leq r \},
    \quad \mbox{   and } \quad
    h = \frac{1}{2}(1,\ldots, 1),\]
  so that
  \[ S = \{ \eps_1 - \eps_2, \eps_2 - \eps_3, \cdots, \eps_{r-1} - \eps_r,
    2\eps_r \} \]
  and $\alpha_r = 2 \eps_r$ is the only simple root not vanishing on~$h$.
  We now observe that
  \[ 2  \eps_r((\eps_{r-1} - \eps_r)^\vee) 
    = 2  \eps_r(\be_{r-1} - \be_r) = - 2 \] 
  is even, so that we do not get any relation, and thus
  $\pi_1(\cO_h)$ is infinite cyclic.\\

  \nin {\bf Non-split type:} This concerns the Lie algebras
  \[ \fsl_n(\H), \quad \fu_{n,n}(\H), \quad \so_{1,d}(\R), d > 2, \quad
    \mbox{ and } \quad \fe_{6(-26)}.\]
  We claim that for all these Lie algebras we have
  $m_{\alpha_j} > 1$. In fact, for $\fsl_n(\bH)$ we have
  $m_\alpha = \dim \H = 4$ for all
  restricted roots. 
  For $\fu_{n,n}(\bH)$ we have $m_{\alpha_j} = 3$
  (type CII in \cite[p.~532]{Hel78}).
  For $\so_{1,d}(\R)$, we have $m_\alpha = d - 1 > 1$, and
  for $\fe_{6(-26)}$ we have $m_\alpha = 8$
    (type EIV in \cite[p.~534]{Hel78}).\\

  \nin {\bf Split type:} This concerns the Lie algebras
  \[ \fsl_n(\R), \quad \so_{p,q}(\R),p,q > 2, \quad
    \fe_6(\R), \quad \mbox{ and }  \quad \fe_7(\R).\]

  For $\fsl_n(\R),  \so_{p,p}(\R), \fe_6(\R)$ and $\fe_7(\R)$,
  these are split simple real Lie algebras of types
  $A_{n-1}, D_p, E_6$ or $E_7$. In particular
  all multiplicities are $1$ and all roots have the same length.
  We therefore have $m_{\alpha_j} = 1$ and there exists an $i$ with 
  $\alpha_j(\alpha_i^\vee) = -1$, so that we obtain the fundamental
  group $\Z_2$.

  For $\so_{p,q}(\R)$ and $p < q$, the situation is more involved.
  If $p + q$ is odd, then the discussion of type BI in
  \cite[p.~532]{Hel78} shows that
  $m_{\alpha_j} = 1$ for $j \leq p-1$ and $m_{\alpha_p} = q - p$.
  As $h = h_1$ in this case, the assumption $p > 2$ entails
  the existence of a simple root $\alpha_2$ with multiplicity~$1$
  and $\alpha_1(\alpha_2^\vee) = -1$.
If $p+q$ is even, this is type DI in \cite[p.~533]{Hel78}
with $m_{\alpha_i} = 1$ for $i < p$ and we can argue as above because~$p > 2$.
\end{rem}

\section{The kernel of the universal complexification} 

A key point in the proof of Theorem~\ref{thm:Z2-struc} is that
the universal complexification $\eta_{\tilde G} \: \tilde G \to \tilde G_\C$ satisfies 
  \begin{equation}
    \label{eq:z2inkereta}
 \tilde Z_2    \subeq \ker(\eta_{\tilde G})
  \end{equation}
follows from  $\eta_{\tilde G}(\tilde Z_2) \subeq 
    \delta_h(\tilde G_\C^h) = \{e\},$ 
  which in turn follows from the connectedness of $\tilde G_\C^h$
  (centralizers of tori in irreducible
  algebraic groups are connected; 
  \cite[Thm.~XIII.4.2]{Ho81}). 
  We now show that we always have equality in \eqref{eq:z2inkereta}.
  
  \begin{thm} \mlabel{thm:4.6} If $G$ is a simple simply connected Lie group
    with Lie algebra $\g$, $h \in \g$ an Euler element,
    and $\eta_G \: G \to G_\C$ its universal complexification, then
    \[ \ker(\eta_G) = Z_2.\]     
  \end{thm}

  \begin{prf} We consider the $4$ types from \cite[Table~3]{MNO23}.

\nin {\bf Complex type:} If $\g$ is of complex type, then
$G$ is a complex group, so that $\eta_G$ is injective.
In this case $Z_2 \subeq \ker(\eta_G) = \{e\}$ (cf.\ \eqref{eq:z2inkereta})
implies equality.

\nin {\bf Non-split type:} If $(\g,h)$ is of non-split type, then
Theorem~\ref{thm:Z2-struc} implies that $Z_2 = \{e\}$.
We therefore have to show that $\eta_G$ is injective,
i.e., that there exists a simply connected matrix Lie group
$G$ with Lie algebra $\g$.

The simple Lie algebras $\g$ for which $(\g,h)$ can be of non-split type are
\[ \fsl_n(\H), \quad \fu_{n,n}(\H), \quad \so_{1,d}(\R), \quad
  \fe_{6(-26)}.\]
For $\g = \fsl_n(\H)$ the corresponding matrix
group $\SL_n(\H)$ with maximal compact subgroup
$\U_n(\H)$ is simply connected. 
The same holds for the matrix group $\U_{n,n}(\H)$ with maximal 
compact subgroup $\U_n(\H) \times \U_n(\H)$ 
(cf.\ also \cite[p.~34]{Ti67}).
For $\g = \so_{1,d}(\R)$, the corresponding simply connected
Lie group $\Spin_{1,d}(\R)$ embeds into the complex spin
group $\Spin_{1+d}(\C)$ (via the spin representation).
For $\g = \fe_{6(-26)}$ we use \cite[p.~46]{Ti67}
to verify that $\eta_G$ is injective. 

\nin {\bf Split type:} If $(\g,h)$ is of split type, then
Theorem~\ref{thm:Z2-struc} implies that $Z_2 \cong \Z_2$.
As $Z_2 \subeq \ker(\eta_G)$ by \eqref{eq:z2inkereta}
and $|\ker(\eta_G)| \leq 2$
by \cite[\S 11.A.2]{Wa92}, the equality $Z_2 = \ker(\eta_G)$ follows.

\nin {\bf Cayley type:} If $(\g,h)$ is of Cayley type,
Theorem~\ref{thm:Z2-struc} implies that $Z_2 \cong \Z$
is generated by $\exp(4 \pi z_1)$, where $z_1 = \frac{i}{2} \gamma_1^\vee$
is a compactly embedded element in an $\fsl_2$-subalgebra corresponding to a
strongly orthogonal root~$\gamma_1$.
As all long roots in $\Sigma_1$ are conjugate under the Weyl group
$\cW_0$ of the root system $\Sigma_0 = \Sigma(\g,\fa) \cap h^\bot$
(the representation of $\g_0$ on $\g_1$ irreducible, so that
the Weyl group acts transitively on set of the extremal weights),
we may assume that $\gamma_1$ is the unique simple root in $\Sigma_1$.

Let $\fk = \g^\theta$, $K = \exp(\fk)$ and recall that
\[ K \cong Z(K)_e \times K' \cong \R \times K'.\]
If $\ft \subeq \fk$ is maximal abelian, then it is a
compactly embedded Cartan subalgebra of $\g$ and
$\Sigma_0$ can be identified with the root system
$\Delta(\fk,\ft) = \Sigma(\fk_\C, i \ft)$. 

For the maximal abelian subgroup $T := \exp_G(\ft) \subeq K$,
and the torus group $T^* := \eta_{G_\C}(T)\subeq G_\C$, we then have
\[ \ker(\exp_T) = 2\pi i \la \Delta(\fk,\ft)^\vee \ra_{\rm grp} 
  \subeq
  \ker(\exp_{T^*}) = 2\pi i \la \Delta(\fg,\ft)^\vee \ra_{\rm grp}\]
(cf.\ \cite[Thm.~12.4.4]{HN12}).
For the unique ``non-compact'' simple root $\gamma_1$, we have
\[ \la \Delta(\fg,\ft)^\vee \ra_{\rm grp}= \la \Delta(\fk,\ft)^\vee \ra_{\rm grp} 
  + \Z \gamma_1^\vee,\]
so that
\begin{equation}
  \label{eq:ker-etag}
\ker(\eta_G)
  \cong \ker(\exp_{T^*})/\ker(\exp_T)
  = \exp(2\pi i \Z \gamma_1^\vee)
  = \exp(4\pi \Z z_1) = Z_2.\qedhere
\end{equation}
This completes the proof.
\end{prf}

\begin{rem} In \cite[\S 11.A.2]{Wa92} it is shown that
  the kernel of the universal complexification
  $\eta_{\tilde G} \: \tilde G \to \tilde G_\C$ is always
  a cyclic group. It is infinite if and only if $\g$ is hermitian,
  and in all other cases it is either trivial or $\Z_2$.
  This information can also be read from the Tits tables \cite{Ti67}, in
  which $\eta_{\tilde G}$ is called $\phi$
  and $\eta_{\tilde G}(\tilde G)$ is called $G^*$, so that
  $\pi_1(G^*) \cong \ker(\eta_{\tilde G})$.
  If $\g$ is not hermitian, then $\ker(\eta_G)$ has at most two elements.
\end{rem}

\end{document}